\documentclass[10pt]{article}

\usepackage{amsmath,amscd,amsthm,makeidx,graphicx,url,bm,bbm,mathrsfs,palatino,bm,mathtools,latexsym}
\usepackage[bbgreekl]{mathbbol}

\usepackage{amssymb}
\usepackage{a4wide,xy,graphicx}
\usepackage{tikz-cd}
\usepackage{enumerate,comment}
\usepackage{titlesec}
\usepackage{tikz-cd}

\numberwithin{equation}{section}

\titleformat{\paragraph}[runin]
{\normalfont\normalsize\bfseries}{\theparagraph}{1em}{}
\titlespacing*{\paragraph}
{0pt}{3.25ex plus 1ex minus .2ex}{1.5ex plus .2ex}

\usepackage[
breaklinks=true,
pdftitle={Very stable Higgs bundles, the nilpotent cone and mirror symmetry},
pdfauthor={Tam\'as Hausel, Nigel Hitchin}]{hyperref}
\usepackage{color}
\definecolor{tocolor}{rgb}{.1,.1,.5}
\definecolor{urlcolor}{rgb}{.2,.2,.6}
\definecolor{linkcolor}{rgb}{.1,.4,.6}
\definecolor{citecolor}{rgb}{.6,.3,.1}
\hypersetup{colorlinks=true, urlcolor=urlcolor, linkcolor=linkcolor, citecolor=citecolor}

\let\oldmarginpar\marginpar
\renewcommand\marginpar[1]{\-\oldmarginpar[\raggedleft\footnotesize #1]
	{\raggedright\footnotesize #1}}

\makeindex

\xyoption{all}
\input{xypic}

\renewcommand{\div}{{\rm div}}

\newcommand{\gen}[1]{\left< #1 \right>} 

\newtheorem{theorem}{Theorem}
\newtheorem{proposition}[theorem]{Proposition}

\newtheorem{lemma}[theorem]{Lemma}

\theoremstyle{remark}
\newtheorem{remark}[theorem]{Remark}

\theoremstyle{definition}
\newtheorem{definition}[theorem]{Definition}

\numberwithin{theorem}{section}

\newcounter{margin}
	{\end{itshape}  \bigskip}

\def\beq{\begin{eqnarray}}
	\def\eeq{\end{eqnarray}}
\def\bes{\begin{eqnarray*}}
	\def\ees{\end{eqnarray*}}

\def\calE{{\mathcal E}}

\DeclareMathOperator{\Aut}{Aut} 

\DeclareMathOperator{\Spec}{Spec} 
\DeclareMathOperator{\Hom}{Hom}

\DeclareMathOperator{\rank}{rank}

\def\C{\mathbb{C}}

\def\calM{{\mathcal{M}}}

\def\M{{\mathcal{M}}}

\def\calA{{\mathcal{A}}}

\def\calF{{\mathcal{F}}}

\def\R{\mathbb{R}}

\newcommand{\T}{{\mathbb{T}}}

\def\t{\mathfrak{t}}

\def\gl{{\mathfrak g\mathfrak l}}

\newcommand{\nc}{\newcommand}

\usepackage{leftidx}

\newcommand{\Gr}{\textnormal{Gr}}

\newcommand{\End}{\textnormal{End}}
\newcommand{\tr}{\textnormal{tr}}

\nc{\op}[1]{\mathop{\mathchoice{\mbox{\rm #1}}{\mbox{\rm #1}}
		{\mbox{\rm \scriptsize #1}}{\mbox{\rm \tiny #1}}}\nolimits}
\nc{\al}{\alpha}

\nc{\ep}{\varepsilon} 
\nc{\ga}{\gamma} 
\nc{\Ga}{\Gamma}
\nc{\la}{\lambda} 
\nc{\La}{\Lambda} 
\nc{\si}{\sigma}
\nc{\Sig}{{\Gamma}} 
\nc{\Om}{\Omega} 
\nc{\om}{\omega}

\nc{\SL}{\mathrm{SL}} 
\nc{\GL}{\mathrm{GL}} 
\nc{\Sp}{\mathrm{Sp}} 
\nc{\SO}{{\mathrm{SO}}}
\nc{\PGL}{\mathrm{PGL}}
\newcommand{\G}{\mathrm{G}}

\nc{\W}{\mathrm{W}}
\nc{\Lg}{\mathrm{L}}
\nc{\Pg}{\mathrm{P}}
\nc{\calL}{{\mathcal L}}
\nc{\Sym}{{\rm Sym}}

\renewcommand{\H}{\mathbb{H}}

\nc{\Frob}{\mathrm{Frob}}
\nc{\spec}{{\rm Spec}}

\def\SU{{\mathrm{SU}}}

\nc{\rN}{\mathrm{N}}

\usepackage{longtable}

\nc{\cpt}{{\op{cpt}}} \nc{\Dol}{{\op{Dol}}} \nc{\DR}{{\op{DR}}}
\nc{\B}{{\op{B}}} \nc{\Triv}{\op{Triv}} \nc{\Hod}{{\op{Hod}}}
\nc{\Log}{{\op{Log}}} \nc{\Exp}{{\op{Exp}}} \nc{\Est}{E_{\op{st}}}
\nc{\Hst}{H_{\op{st}}} \nc{\Left}[1]{\hbox{$\left#1\vbox to
		10.5pt{}\right.\nulldelimiterspace=0pt \mathsurround=0pt$}}
\nc{\Right}[1]{\hbox{$\left.\vbox to
		10.5pt{}\right#1\nulldelimiterspace=0pt \mathsurround=0pt$}}
\nc{\LEFT}[1]{\hbox{$\left#1\vbox to
		15.5pt{}\right.\nulldelimiterspace=0pt \mathsurround=0pt$}}
\nc{\RIGHT}[1]{\hbox{$\left.\vbox to
		15.5pt{}\right#1\nulldelimiterspace=0pt \mathsurround=0pt$}}

\nc{\bee}{{\bf E}} 

\title{Hitchin map on even very stable upward flows}
\author{Miguel Gonz\'alez \\ {\it ICMAT Madrid} 
\\{\tt mgonzalez.contacto@gmail.com} \and  Tam\'as Hausel \\ {\it IST Austria} 
\\{\tt tamas.hausel@ist.ac.at}}

\begin{document}

\maketitle

\begin{abstract}We define even very stable Higgs bundles and study the Hitchin map restricted to their upward flows. In the $\GL_n$ case we classify the type $(1,\dots,1)$ examples, and find that they are governed by a root system formed by the roots of even height. We  discuss how the spectrum of equivariant cohomology of real and quaternionic Grassmannians, $4n$-spheres and the real Cayley plane appear to describe the Hitchin map on even cominuscule upward flows. The even upward flows in question are the same as upward flows in Higgs bundle moduli spaces for quasi-split inner real forms. The latter spaces have been pioneered by Oscar García-Prada and his collaborators.
\end{abstract}

\section{Introduction}

This paper is a write-up of the second author's talk \cite{Hausel22} at the conference "Moduli spaces and geometric structures"
 in honour of Oscar García-Prada on the occasion of his 60th birthday at ICMAT Madrid in September 2022.
 
In \cite{HauselHitchin}, motivated by mirror symmetry, the notion of very stable Higgs bundle was introduced. Let $C$ be a smooth projective curve. Let $\M$ denote the moduli space of rank $n$ degree $d$ semistable Higgs bundles $(E,\Phi)$, where $E$ is a rank $n$ degree vector bundle and $\Phi\in H^0(C;\End(E)\otimes K)$ is a Higgs field. There is a $\T$-action on $\M$ by scaling the Higgs field, i.e. $\lambda \in \T$ acts by sending $(E,\Phi)$ to $(E,\lambda\Phi)$. A fixed point  $\calE\in \M^{s\T}$ is called {\em very stable}, if the upward flow
$$W^+_\calE:=\{\calF\in \M : \lim_{\lambda\to 0} \lambda\cdot\calF=\calE\} \subset \M$$ is closed. In Section~\ref{bb} we recall the basic properties of very stable upward flows in general as well as for the the moduli space of Higgs bundles $\M$. 

One of the main results of \cite{HauselHitchin} is the classification of very stable Higgs bundles $(E,\Phi)\in \M^{\T}$ of type $(1,\dots,1)$. A fixed point is of type $(1,\dots,1)$ when the vector bundle $E=L_0\oplus\dots \oplus L_{n-1}$ is a direct sum of line bundles, and the Higgs field  $\Phi|_{L_i}:L_i \to L_{i+1}K\subset E K$, which we denote by \begin{align*}  \label{bi} b_i:=\Phi|_{L_i}\in \Hom(L_i,L_{i+1}K)\cong H^0(C;L_i^{-1}L_{i+1}K).\end{align*}  Then we have
\begin{theorem}[{\cite[Theorem 4.16]{HauselHitchin}}] The type $(1,...,1)$ Higgs bundle $(E,\Phi)\in \M^\T$ is very stable if and only if the divisor $\div(b_1)+\dots+\div(b_{n-1})$ is reduced.
\end{theorem} We recall this classification in Theorem~\ref{classification} below, and a reformulation of it in Remark~\ref{minusculeremark} in terms of minuscule dominant weights of $\GL(n,\C)$.

Garcia-Prada and Ramanan in \cite{garciaprada_involutions_2019} study involutions on the moduli space of Higgs bundles. One important involution $\theta:\M\to \M$ is given by $\theta(E,\Phi):=(E,-\Phi)$. In \cite{garciaprada_involutions_2019} it is shown that the fixed points $\M^\theta$ correspond to $U(p,n-p)$-Higgs bundles (including the case $p=0$, where $U(0,n) := U(n)$). We recall these notions in Section~\ref{even}. 

In this paper we will be interested in the so-called {\em even  upward flows} $W^{2+}_\calE$ for any $\calE\in \M^{s\T}\subset \M^\theta$ which are defined to be the upward flows $\calE$ in the semi-projective $\M^\theta$, or equivalently, the intersection  $W^{2+}_\calE:=W^+_\calE\cap \calM^\theta$. Then we can define {\em even very stable} Higgs bundles $\calE\in \M^{\T}$ for which the even upward flow $W^{2+}_\calE\subset \M^\theta$ is closed. One of the main results of this paper is the following

\begin{theorem}  The type $(1,...,1)$ Higgs bundle $(E,\Phi)\in \M^\T$ is even very stable, if and only if the divisors $\div(b_2)+\dots+\div(b_{n-2})$ and  $\div(b_i)+\div(b_{i+2k+1})$ for
$1\leq i \leq i+2k+1\leq n-1$ are all reduced. 
\end{theorem}

To clarify the meaning of this complicated looking set of divisors, we reformulate this theorem in Theorem~\ref{evenmin} in terms of so-called {\em even minuscule} dominant weights using positive weights of even height.

As the Hitchin map restricted to very stable upward flows is finite flat and $\T$-equivariant between affine spaces, with positive $\T$-action of the same dimension it is suspectible of explicit description. In the type $(1,\dots,1)$ very stable case, the second author found  such an explicit description in \cite{Hausel-hmsec} in terms of the spectrum of equivariant cohomology of the Grassmannian $\Gr_k(\C^n)$. We will recall this in Section~\ref{modell} below. 

Finally in Section~\ref{thetamodell} we study the problem of modelling the  Hitchin map on certain even very stable upward flows, in terms of the equivariant cohomology of  homogeneous spaces. We will find in Theorem~\ref{final} that for $\GL_{2n}$ the equivariant cohomology of quaternionic Grassmannians, for $\SO_{4n+2}$ the equivariant cohomology of the $4n$-sphere and finally for ${\rm E}_6$ the equivariant cohomology of the real Cayley plane should model the Hitchin map on some specific even very stable flows. The appearance of these symmetric spaces is interesting, partly because they are not of Hermitian type, and also because they are quotients of the Nadler group \cite{Nadler} of the quasi-split real form of Hodge type  (see \cite[Section 2.3]{garciaprada_involutions_2019} for the definition). 

In this paper we are concentrating on type $(1,...,1)$ very stable and even very stable upward flows. By now there are many interesting results about other types of very stable or wobbly Higgs bundles see e.g. \cite{HitchinMultiplicity} for multiplicity algebras of type $(2)$ very stable Higgs bundles, \cite{Franco-etal} for many wobbly Higgs bundles - both papers in this conference proceedings - and   \cite{Peon-Nieto} for a classification of all type $(n_1,n_2)$ very stable components. 

 {\noindent \bf Acknowledgements.} Most of the research for this paper was done when the first author visited the second author's group at IST Austria as a summer intern in 2022. The first author is grateful for the hospitality and support received during this stay. We thank Oscar Garc\'ia-Prada for engineering the internship, and for constant support. We also thank
 Andreas \v{C}ap, Mischa Elkner, Tim Henke, Nigel Hitchin, Friedrich Knop, Jakub L\"owit,  David Nadler, Ana Pe\'on-Nieto, Kamil Rychlewicz and Anna Sis\'ak for useful discussions.  The second author was supported by an FWF grant ``Geometry of the top of the nilpotent cone'' number P 35847.

\section{Bialynicki-Birula decomposition}\label{bb}

In this section we first recall the definition of a semi-projective variety and then collect the basics of the Bialynicki-Birula decomposition associated to such a variety. 

\begin{definition}
	
	Let $X$ be a normal complex quasi-projective variety equipped with a $\mathbb T:=\mathbb C^\times$ action. $X$ is {\em semi-projective} if the fixed point locus $X^{\mathbb T}$ is projective, and for every $x \in X$ there is a $p \in X^{\mathbb T}$ such that $\lim\limits_{\lambda \to 0}\lambda x = p$.
\end{definition}

The latter is to be understood as the existence of a $\mathbb T$-equivariant morphism $f: \mathbb A^1 \to X$ such that $f(1) = x$ and $f(0) = p$. Semi-projective varieties are endowed with a stratification in affine subvarieties known as the {\em Bialynicki-Birula decomposition} \cite{BialynickiBirula}, which we now recall. We refer to \cite[Section 2]{HauselHitchin} for further details.

\begin{definition}
	Let $X$ be a semi-projective variety and $\alpha \in X^{\mathbb T}$. The {\em upward flow from $\alpha$} is defined to be
	$$W_{\alpha}^+ := \{x \in X : \lim\limits_{\lambda \to 0}\lambda x = \alpha\}.$$
	
	Similarly, the {\em downward flow from $\alpha$} is 
	$$W_{\alpha}^- := \{x \in X : \lim\limits_{\lambda \to \infty}\lambda x = \alpha\}.$$
	
	The {\em Bialynicki-Birula partition} is $X= \bigsqcup_{\alpha \in X^{\mathbb T}}W_\alpha^+$. The {\em core} of $X$ is defined to be $\mathcal C:=\bigsqcup_{\alpha \in X^{\mathbb T}}W_\alpha^-$.
	
\end{definition}
 
 \begin{definition}
	For a connected component of the fixed locus, $F \in \pi_0(X^{\mathbb T})$, we define its {\em attractor} as $W^+_F := \bigcup_{\alpha \in F}W_\alpha^+$, and its {\em repeller} as $W^-_F := \bigcup_{\alpha \in F}W_\alpha^-$. The {\em Bialynicki-Birula decomposition} is $X= \bigsqcup_{F \in \pi_0(X^{\mathbb T})}W_F^+$.
 \end{definition}

 Given a smooth fixed point $\alpha \in X^{s\mathbb T}$, the $\mathbb T$-action on $X$ induces a representation of $\mathbb T$ on the tangent space $T_\alpha X$. We denote, for $k \in \mathbb Z$, the weight space $(T_\alpha X)_k \subset T_\alpha X$ where $\lambda \in \mathbb T$ acts via multiplication by $\lambda^k$. This leads to a decomposition $T_\alpha X = \bigoplus_{k \in \mathbb Z} (T_\alpha X)_k$ in weight spaces. We denote $T_\alpha^+X := \bigoplus_{k>0}(T_\alpha X)_k$ the positive part and $T_\alpha^-X := \bigoplus_{k<0}(T_\alpha X)_k$ the negative part. We have:
 
 \begin{proposition}
 	Given a smooth fixed point $\alpha \in X^{s \mathbb T}$, the upward flow $W_\alpha^+$ (resp. the downward flow $W_\alpha^-$) is a locally closed $\mathbb T$-invariant subvariety of $X$ which is isomorphic to $T^+_\alpha X$ (resp. $T^-_\alpha X$) as varieties with $\mathbb T$-action.
 \end{proposition}

The proof was originally given in \cite{BialynickiBirula} for smooth complete $X$. A proof for the general case is given in \cite[Proposition 2.1]{HauselHitchin}. 

Finally, suppose further that $X^s$ is equipped with a symplectic form $\omega \in \Omega^2(X^s)$ such that, for $\lambda \in \mathbb T$, we have $\lambda^*(\omega) = \lambda\omega$. This supposition is motivated by the fact that the semi-projective variety we will be studying, the moduli space of semistable Higgs bundles, is endowed with such a form. Then, we have:

\begin{proposition}
	For a smooth point $\alpha \in X^{s\mathbb T}$, the subspaces $T^+_\alpha X$ and $T^{\le 0}_\alpha X := (T_\alpha X)_0 \oplus T^-_\alpha X$ of $T_\alpha X$ are Lagrangian. Moreover, the subvarieties $W_\alpha^+$ and $W_{F_\alpha}^-$ are also Lagrangian.
\end{proposition}   

	The proof is given in \cite[Proposition 2.10]{HauselHitchin}. The main idea is that, for $v \in (T_\alpha X)_k$ and $w \in (T_\alpha X)_l$, we have 
	$$\lambda\omega(v,w) = \lambda^*(\omega)(v,w) = \omega(\lambda \cdot v, \lambda \cdot w) = \omega(\lambda^kv, \lambda^lw) = \lambda^{k+l}\omega(v,w),$$
	
	\noindent so that $\omega(v,w)$ can only be nonzero in the situation $k+l=1$, which does not happen if $k,l>0$ or if $k,l\le 0$. 

	\begin{definition}
		We say that $\alpha \in X^{s\mathbb T}$ is {\em very stable} if $W^+_\alpha \cap \mathcal C = \{\alpha\}$.
	\end{definition}

	This definition was introduced in \cite[Definition 4.1]{HauselHitchin}, where it was proven \cite[Lemma 4.4]{HauselHitchin} that $\alpha \in X^{s\mathbb T}$ is very stable if and only if $W^+_\alpha \subset X$ is closed.

\subsection{Lagrangian upward flows in $\M$}

	In this section we introduce Higgs bundles and show how the previous theory of Bialynicki-Birula applies to the moduli space of semistable Higgs bundles. For this, we fix a smooth projective curve $C$ over the complex numbers with genus $g \ge 2$ and canonical line bundle $K$. 
	
	\begin{definition}
		A {\em Higgs bundle} is a pair $(E,\Phi)$ where $E$ is a holomorphic vector bundle over $C$ and $\Phi \in H^0(C, \End(E) \otimes K)$.
	\end{definition} 

	Such an object can be defined in more generality for a real reductive Lie group $G$ \cite[Definition 3.1]{HiggsPairs}, giving $G$-Higgs bundles.
	The above definition is recovered by setting $G = GL_n(\mathbb C)$ for $n = \rank E$. Recall that a Higgs bundle $(E,\Phi)$ is {\em stable} if, for every nonzero proper vector subbundle $F \subset E$ such that $\Phi(F) \subseteq F \otimes K$, we have
	$$\mu(F):=\frac{\deg F}{\rank F} < \mu (E) := \frac{\deg E}{\rank E},$$

	\noindent and it is {\em semistable} if for the same subbundles we have $\mu(F) \le \mu(E)$. We denote by $\mathcal M := \mathcal M^d_n$ the moduli space of semistable Higgs bundles of fixed rank $n$ and degree $d$. It was constructed via gauge theory in \cite{HitchinSelfDuality} and by algebraic geometric methods in \cite{Nitsure,Simpson}. It is a normal \cite{SimpsonII} quasi-projective variety with a hyperkähler metric at its smooth points, which are the stable Higgs bundles. In particular, it has a symplectic structure $\omega \in \Omega^2( \mathcal M^s)$.
	
	This space also carries a natural $\mathbb T$-action defined by $(E,\Phi) \mapsto (E,\lambda \Phi)$ which turns it into a semi-projective variety and such that $\lambda^*(\omega) = \lambda\omega$. Thus, the Bialynicki-Birula theory from the previous section applies. The fixed locus $\mathcal M^{\mathbb T}$ can be identified as follows. We have, for any $\lambda \in \mathbb T$, an isomorphism of vector bundles $f_\lambda \in \Aut(E)$ such that
	\begin{equation}\label{lambdaphi}
		\Phi \circ f_\lambda = f_\lambda \circ (\lambda \Phi).
	\end{equation}

	In other words, we have an action of $\mathbb T$ on $E$ which is linear on each fiber. Hence we can decompose $E = L_0 \oplus \dots \oplus L_k$ into weight spaces, where $f_\lambda|_{L_i} = \lambda^{w_i} \cdot \text{Id}_{L_i}$. The compatibility condition (\ref{lambdaphi}) shows that, if $v_i \in L_i$, then $f_\lambda(\lambda \cdot \Phi(v_i)) = \Phi(f_\lambda(v_i)) = \lambda^{w_i}\Phi(v_i)$, thus $f_\lambda(\Phi(v_i)) = \lambda^{w_i-1}\Phi(v_i)$. Hence, $\Phi$ maps the space for weight $w_i$ into the space for weight $w_i-1$. In particular, the weights can be chosen to be of the form $w_i = w_0 - i$ and the Higgs field has the property $\Phi(L_i) \subseteq L_{i+1}\otimes K$. From this, we can associate an invariant $(\rank L_0, \dots, \rank L_k)$ to the fixed point, known as the {\em type}.
	
	The previous decomposition shows that, in fact, a Higgs bundle fixed by the $\mathbb T$-action is nilpotent, since $\Phi^{k+1} \equiv 0$. Another way of seeing this is via the {\em Hitchin map}:
	$$h: \mathcal M \to \mathcal A:=\bigoplus_{i=1}^nH^0(C, K^i)$$
	
	\noindent defined by the coefficients $a_i \in H^0(C, K^i)$ of the characteristic polynomial $\det(\Phi - xI_n) = x^n + \sum_{j=1}^{n}a_{j}x^{n-j}$. This map is a proper, completely integrable Hamiltonian system \cite{HitchinIntegrable, Nitsure} whose fibers are Lagrangian at their smooth points, and the generic fibers are abelian varieties \cite{BNR}. Moreover, by letting $\mathbb T$ act on $H^0(C,K^i)$ with weight $i$, the Hitchin map is $\mathbb T$-equivariant. Thus, $\mathcal M^{\mathbb T} \subseteq h^{-1}(0)$ so that fixed points $(E,\Phi)$ have characteristic polynomial $x^n$ and are nilpotent.
	
	The upward and downward flows from the Bialynicki-Birula partition have been characterized in \cite[Proposition 3.4 and 3.11]{HauselHitchin}:
	
	\begin{proposition}\label{upwardflow}
		Let $\mathcal E = (E',\Phi') \in \mathcal M^{s\mathbb T}$ and $(E,\Phi) \in \mathcal M$. We have $(E',\Phi') \in W_{\mathcal E}^+$ if and only if there exists a filtration
		$$0 = E_0 \subset E_1 \subset \dots \subset E_k = E$$
		
		\noindent such that $\Phi(E_i) \subseteq E_{i+1} \otimes K$ and the associated graded object verifies $(\Gr{E}, \Gr(\Phi)) \simeq (E',\Phi')$. The same is true replacing $W^+_{\mathcal E}$ with $W^-_{\mathcal E}$ and the ascending filtration with a descending filtration.
	\end{proposition}
	
	The downward flows have a concrete characterization via the Hitchin map. Since the action of $\mathbb T$ on $\mathcal A$ is by positive weights, the core is just $\mathcal C_{\mathcal A} = \{0\}$. Because $h$ is $\mathbb T$-equivariant, this shows that $\mathcal C := \mathcal C_{\mathcal M} \subset h^{-1}(0)$. On the other hand, the properness of $h$ shows that $h^{-1}(0)$ is projective and thus $\mathcal C = h^{-1}(0)$. This is typically called the {\em nilpotent cone} in this context. Notice that $\T$-equivariance of $h$ implies $\mathcal M^{\T} \subseteq \mathcal C$. From this, the notion of being very stable becomes:
	
	\begin{definition}
		A Higgs bundle $\mathcal E = (E,\Phi) \in \mathcal M^{s\mathbb T}$ is {\em very stable} if the only nilpotent Higgs bundle in $W_{\mathcal E}^+$ is $\mathcal E$ itself. Otherwise, it is {\em wobbly}.
	\end{definition}
	
	One interesting aspect of such objects is that the Hitchin map restricts nicely to the upward flow \cite[Lemma 4.6]{HauselHitchin}:
	
	\begin{proposition}
		If $\mathcal E \in \mathcal M^{s\mathbb T}$ is very stable then $h: W^+_\mathcal E \to \mathcal A$ is finite, flat, surjective and generically étale.
	\end{proposition}

	
	\subsection{Examples of very stable Higgs bundles}\label{examples}
	
	In this section we recall from \cite{HauselHitchin} some examples of very stable Higgs bundles. First we consider the fixed point component of type $(n)$. The fixed points of this type are elements of the form $(E,0)$ with $E$ a semistable Higgs bundle. Thus, this component is just the moduli space of semistable rank $n$ degree $d$ vector bundles, $\mathcal N$. The upward flow for $\mathcal E = (E,0)$ is given by $W_{\mathcal E}^+ = \{(E,\Phi) : \Phi \in H^0(C, \End(E) \otimes K)\}$, so that $\mathcal E$ is very stable if and only if the only nilpotent Higgs field $\Phi \in H^0(C, \End(E) \otimes K)$ it admits is $\Phi \equiv 0$. This is the notion of very stable vector bundle introduced by Drinfeld and Laumon \cite{Laumon}, for which they prove that very stable bundles form an open dense subset of the component. 
	
	Next, we shall focus on the type $(1,1,\dots,1)$ case. The starting example of very stable Higgs bundle in this component is the {\em canonical uniformising Higgs bundle}, $\mathcal E_0 = (E_0,\Phi_0)$, where
	$$E_0 = \mathcal O \oplus K^{-1} \oplus \dots \oplus K^{1-n},$$
	
	\noindent and, given $a = (a_1,\dots,a_n) \in \mathcal A = H^0(C,K) \oplus \dots \oplus H^0(C,K^n)$, the Higgs field
	$$\Phi_a = \begin{pmatrix}
		0 & 0 & \dots & 0 & a_n \\
		1 & 0 & \dots & 0 & a_{n-1} \\
		0 & 1 & \dots & 0 & a_{n-2} \\
		\vdots & \vdots & \ddots & \vdots & \vdots \\
		0 & 0 & \dots & 1 & a_1\end{pmatrix}$$

	\noindent is given by the companion matrix. The map $a \mapsto (E_0, \Phi_a)$ provides a section of the Hitchin map, known as the {\em Hitchin section} \cite{HitchinSection}. By means of Proposition \ref{upwardflow}, it follows that $\{(E_0,\Phi_a) : a \in \mathcal A\} \subseteq W^+_{{\mathcal E}_0}$. Moreover, since both are affine spaces of equal dimension $\dim\mathcal M/2$, the upward flow is precisely the Hitchin section, hence ${\mathcal E}_0$ is very stable.   
	
	It is possible to completely classify very stable Higgs bundles of this type by starting with this example and performing Hecke transformations. First we note that the data of a type $(1,1,\dots,1)$ fixed point is equivalent to the choice of a line bundle $L_0$ over $C$ (that is, a divisor $\delta_0$ up to principal divisor), as well as effective divisors $\delta_1,\dots,\delta_{n-1}$ on $C$. Indeed, such a fixed point $(E,\Phi)$ is of the form $E = L_0 \oplus \dots \oplus L_{n-1}$, where for all $j$ we have $\rank L_j = 1$, and $\Phi|_{L_{j-1}} = b_{j-1}$ for nonconstant maps $b_j : L_{j-1}\to L_{j} \otimes K$, $1 \le j \le n-1$. Thus, $L_0$ is given, and $\delta_i$ is obtained as the zero locus of $b_i$ with multiplicities. On the other hand, given $(\delta_0,\delta_1,\dots,\delta_{n-1})$ we construct $E$ by setting $L_i := L_0 \otimes \mathcal O(\delta_1+\dots+\delta_{i-1}) \otimes K^{-i}$ and $b_i := s_{\delta_i} \in \mathcal O(\delta_i) = \mathcal O(L_{j-1}^*\otimes L_j \otimes K)$ the canonical section. We shall denote the bundle corresponding to $\delta := (\delta_0,\dots,\delta_{n-1})$ by $\mathcal E_{\delta} = (E_\delta, \Phi_\delta)$. 
	
	Another convenient way of labelling these points is via choosing a dominant weight of $\GL(n,\mathbb C)$ at each point of $C$, that is, a map 
	$$\mu : C \to \Lambda^+(\GL(n, \mathbb C)) = \left\{\sum_{i=1}^{n}a_i\omega_i : a \in \mathbb Z_{\ge 0}^{n-1} \times \mathbb Z\right\},$$ 
	
	\noindent where the $\omega_i$ are the fundamental weights. We require that the set $C \setminus \{\mu = 0\}$ is finite. We then define 
	$$\delta_\mu :=\left(\sum_{c\in C}\langle {\mu}(c),\omega_{n}^\vee \rangle c ,\left(\sum_{c\in C} \langle {\mu}(c),\omega_{i}^\vee \rangle c \right)_{i=1,\dots,n-1}\right),$$ 
	
	\noindent and $\mathcal E_\mu := \mathcal E_{\delta_\mu}$. Conversely, we can retrieve the map $\mu$ from $\delta$ as $${\mu}_{\delta}(c)=\delta_0(c)\omega_n+\sum_{i=1}^{n-1}\delta_{i}(c)\omega_{i}\in \Lambda^+,$$ where $D(c)$ for a divisor $D$ means the coefficient of $c$ in $D$.
	
	\subsection{Hecke transformations}
	Now we explain Hecke transformations for Higgs bundles. These play a key role since they allow to relate the upward flows of the different fixed points in the type $(1,1,\dots,1)$ component. In order to define Hecke transformations of $(E,\Phi)$, we start by choosing a point $c \in C$ and a subspace $V \in \Gr(k, E|_c)$ which is $\Phi|_c$-invariant, that is, $\Phi|_c(V) \subset V \otimes K|_c$. The Hecke transformation $\mathcal H_V(E,\Phi) := (E',\Phi')$ is defined by diagram
	$$\begin{array}{ccccccc} 0 \to & E^\prime& {\to} &E& {\to} & E|_c/V & \to 0 
		\\ & \Phi^\prime \downarrow &&\Phi \downarrow&& \overline{\Phi}_c\downarrow & \\ 0 \to & E^\prime \otimes  K&\to &E \otimes K&\to& E|_c/V \otimes K&\to 0 \end{array},$$
	\noindent where $E|_c/V$ is to be regarded as a skyscraper sheaf at $c$. More details of the construction of this diagram can be found in \cite[Definition 4.10]{HauselHitchin}.
	
	It is possible to reach any $\mathcal E_\mu$ from successive Hecke transformations that start at $\mathcal E_0$. This is due to the following fundamental operation: starting with $\mathcal E_{0}$ and selecting the natural invariant subspace $V_k = (L_k \oplus \dots \oplus L_{n-1})|_c$, the resulting Hecke transformation gives $\mathcal E_{\mu_{c,k}}$, where $\mu_{c,k}(c) = \omega_k$ and zero otherwise. This is explained in \cite[Example 4.13]{HauselHitchin}. For arbitrary $\mu$ it suffices to iterate the previous operation for every $c \in C$, at the $\omega_k$ indicated by $\mu(c)$. One of the main results of $\cite{HauselHitchin}$ is that the upward flows are also related by Hecke transformations, from which the following classification can be deduced:
	
	\begin{theorem}[{\cite[Theorem 4.16]{HauselHitchin}}]\label{classification}
		A stable fixed point of type $(1,1,\dots,1)$, $(E_\delta, \Phi_\delta) \in \mathcal M^{s\mathbb T}$, is very stable if and only if the divisor $\delta_1 + \dots + \delta_{n-1}$ is reduced.
	\end{theorem}

	\begin{remark}\label{minusculeremark}
		The previous statement can be rephrased as $(E_\mu, \Phi_\mu)$ being very stable if and only if for every $c \in C$, either $\mu(c) = a_n\omega_n$ or $\mu(c) = a_n\omega_n + \omega_k$ where $k \in \{1,\dots,n-1\}$ and $a_n \in \mathbb Z$. In other words, the point is very stable if and only if for every $c \in C$, the weight $\mu(c)$ is {\em minuscule}, that is, minimal with respect to the partial ordering in $\Lambda^+(\GL(n, \mathbb C))$  given by $\mu_1 \ge \mu_2 \iff \mu_1-\mu_2 \in \Phi^+ = \{\sum_{\alpha \in \Delta^+}a_\alpha \alpha : a_\alpha \in \mathbb Z_{\ge 0}\}$, where $\Delta^+$ denotes the set of positive roots. 
	\end{remark}

\section{Even very stable upward flows} \label{even}


	In this section we extend the results summarized above to the subspace $\mathcal M^\theta \subset \mathcal M$ of the moduli space defined by the fixed points of the subgroup $C_2 = \{1,-1\} \subseteq \mathbb T$, acting as the involution $\theta : (E,\Phi) \mapsto (E,-\Phi)$. Clearly, this subspace contains all the fixed points by the $\mathbb T$-action, so the previous concepts can be extended naturally. Moreover, by \cite[Theorem 6.3]{garciaprada_involutions_2019}, the space $\mathcal M^\theta$ contains the images of the maps $\mathcal M_{U(p,q)} \rightarrow \mathcal M$, given by extension of structure group, of the moduli spaces of $U(p,q)$-Higgs bundles (for the different $U(p,q)$ with $p+q=n$, $p\le q$) into the moduli space of $\GL(n,\mathbb C)$-Higgs bundles. The stable locus $\mathcal M^{\theta,s}$ is covered by these images, which are the Higgs bundles that we will consider.
	
	We start by recalling the following definition from \cite[Definition 3.3]{bradlow_surface_2003}.

	\begin{definition}
		A {\em $U(p,q)$-Higgs bundle} $(E, \Phi)$ is a holomorphic vector bundle $E$ of the form $E = V \oplus W$, where $V$ and $W$ are vector bundles of ranks $p$ and $q$, respectively, and $\Phi \in H^0(\End(E)\otimes K)$ is a section satisfying $\Phi(V) \subset W \otimes K$, $\Phi(W) \subset V \otimes K$.
	\end{definition}

    We denote by $\mathcal M^s_{U(p,q)}$ the moduli space of stable $U(p,q)$-Higgs bundles, where stability is defined as for $\GL(n,\mathbb C)$-Higgs bundles. As explained before, this space sits inside $\mathcal M$ as the fixed point locus of the involution $\theta$, that is, a stable Higgs bundle $(E,\Phi)$ is a $U(p,q)$-Higgs bundle (for some $p$ and $q$) if and only if $(E,\Phi) \simeq (E,-\Phi)$ \cite[Theorem 6.3]{garciaprada_involutions_2019}. We will often use interchangeably the moduli space $\mathcal M_{U(p,q)}$ and its image inside $\mathcal M$.
    
    Note that every fixed point $\mathcal E = (E, \Phi)$ of the $\mathbb T$-action is in particular fixed by $-1$ and hence a $U(p,q)$-Higgs bundle for some $p$ and $q$. This can be seen more explicitly by observing that, in the decomposition $E = E_0 \oplus \dots \oplus E_{k-1}$ with $\Phi(E_i) \subset E_{i+1} \otimes K$, the Higgs field $\Phi$ interchanges the summands with odd indices by those of even indices. In particular, type $(1,1,\dots,1)$ fixed points are Higgs bundles for the quasi-split group $U(p,p)$ or $U(p,p+1)$. 

    Moreover, recall from Białynicki-Birula theory in Section \ref{bb}
    that the upward flow $W^+_{\mathcal E}$ is an affine space isomorphic to the subspace of positive weights $T^+_\mathcal E\mathcal M \subset T_\mathcal E\mathcal M$. It is easy to identify the subspace of $U(p,q)$-Higgs bundles.
	
	\begin{remark}
		The $U(p,q)$-Higgs bundles in $W^+_{\mathcal E}$, that is, $\mathcal M_{U(p,q)} \cap W^+_{\mathcal E}$, correspond via the isomorphism $W^+_{\mathcal E} \simeq T^+_\mathcal E\mathcal M$ to the vector subspace of positive, {\em even} weights $T^{2+}_{\mathcal E}\mathcal M \subset T_{\mathcal E}\mathcal M$.
	\end{remark}
 
 	\noindent This is because the even weights are precisely the vector subspace fixed by multiplication by $-1$, and the previous isomorphism is $\mathbb T$-equivariant. In other words, we can view the locus of $U(p,q)$-Higgs bundles at the upward flow of a fixed point as a subspace:
	
	\begin{definition}
		The {\em even upward flow} at $\mathcal E$ is the subspace of $W^+_{\mathcal E}$ corresponding to $T^{2+}_{\mathcal E}\mathcal M$ and is denoted by $W^{2+}_{\mathcal E}$.
	\end{definition}

	Note that this coincides with the standard upward flow when defined in $\mathcal M^\theta$ instead of $\mathcal M$. Hence, we also have the following natural definition of very stable points: 

	\begin{definition}
		We say that $\mathcal E$ is {\em even very stable} if $W^{2+}_{\mathcal E} \cap \mathcal C = \{\mathcal E\}$, where $\mathcal C = h^{-1}(0) \subseteq \mathcal M$ denotes the locus of nilpotent Higgs bundles. Otherwise, it is said to be \textit{even wobbly}.
	\end{definition}

	We remark, as one of the main interests for this study, that the subspaces of even weights $T^{2}_{\mathcal E}\mathcal M \subset T_{\mathcal E}\mathcal M$ that we are considering are Lagrangian, since the symplectic form $\omega$ pairs the subspace of weight $k$ with that of $1-k$, as explained at the end of Section \ref{bb}, so that the subspaces of even weights are paired with those of odd weights. In fact, the subvariety $\mathcal M^\theta \subset \mathcal M$ itself is Lagrangian, as explained in \cite[Theorem 8.10]{garciaprada_involutions_2019}.
	
	Obviously, a very stable fixed point is also even very stable. However, a wobbly fixed point $\mathcal E$ can either remain even wobbly or instead be even very stable, depending on whether the nontrivial intersection $(W^+_{\mathcal E} \cap \mathcal C) \setminus \{\mathcal E\}$ happens at even weights or not. We will classify the even very stable Higgs bundles of type $(1,1,\dots,1)$, revealing that both situations already arise in this case. We follow the same notation of Subsection \ref{examples}.

 \begin{proposition}\label{wobbly1}
		Let $\mathcal E = (E, \Phi)$ be a smooth fixed point of type $(1,1,\dots,1)$
        . Suppose that there is some point $c \in C$ such that $b_i(c) = b_j(c) = 0$ for some $i, j$ of different parity. Then $\mathcal E$ is even wobbly.
	\end{proposition}

	\noindent {\it {\it Proof.}} The proof follows the approach of constructing a curve via Hecke transformations followed in the proof of \cite[Theorem 4.16]{HauselHitchin}. We will see that, in this situation, the curve can be constructed along even weights. We have that $b_k(c) = b_{k+l}(c) = 0$ where $l > 0$ and odd. We start by performing a Hecke transformation at the $\Phi_c$-invariant subspace $(L_0 \oplus \dots \oplus L_{k-1})|_c \subset E|_c$, yielding 
	
	$$E' = L_0 \oplus \dots \oplus L_{k-1} \oplus L_k(-c) \oplus \dots \oplus L_{n-1}(-c),$$
	
	$$\Phi' = \begin{pmatrix}
	0 & 0 & \dots & 0 & \dots &0& 0 \\
	b_1 & 0 & \dots & 0 & \dots &0& 0 \\
	0 & b_2 & \dots & 0 & \dots &0& 0 \\
	\vdots & \vdots & \ddots & \vdots & \ddots & \vdots & \vdots \\
	0 & 0 & \dots & \frac{b_k}{s_c} & \dots & 0& 0\\
	\vdots & \vdots & \ddots & \vdots & \ddots & \vdots & \vdots \\
	0 & 0 & \dots & 0 & \dots & b_{n-1} & 0
	\end{pmatrix},$$
	
	\noindent where $s_c \in H^0(\mathcal O(c))$ is the canonical section. This bundle is still stable \cite[Lemma 4.17]{HauselHitchin}. 
    We have that $V_0 = (L_{k}(-c) \oplus \dots \oplus L_{n-1}(-c))|_c \subset E'|_c$ is the $(n-k)$-dimensional $\Phi'_c$-invariant subspace that transforms $(E', \Phi')$ back into $(E(-c), \Phi)$. Now let $\{v_0, \dots, v_{n-1}\}$ be a basis of $E'|_c$, each $v_j$ taken in the corresponding component. We define the following $(n-k)$-dimensional subspace:
	$$V := \gen{v_{k-1} + v_{k+l}, v_k, v_{k+1}, \dots,v_{k+l-1},v_{k+l+1}, v_{n-1}} \subset E'|_c$$
	
	This is another $\Phi'_c$-invariant subspace, since $\Phi'_c(v_{k-1}+v_{k+l}) \in \gen{v_k, v_{k+l+1}}$, $\Phi'_c(v_m) \in \gen{v_{m+1}}$ for $m \in \{k, \dots, k+l-2\} \cup \{k+l+1, \dots, n-1\}$, and $\Phi'_c(v_{k+l-1}) = 0$ (in this analysis we take $v_j = 0$ for any $j > n-1$). Now, recall that there is an induced $\mathbb T$-action on $Gr_{n-k}(E'|_c)$ given by the $\mathbb T$-action on $E'$ with weight $-i$ on the $i$-th summand of $E'$. With this, given $\lambda \in \mathbb T$ we define $V_\lambda := \lambda V$, that is:
	$$V_\lambda = \gen{\lambda^{l+1}v_{k-1} + v_{k+l}, v_k, v_{k+1}, \dots,v_{k+l-1},v_{k+l+1}, v_{n-1}}.$$
	
	Note that this subspace is always $n-k$ dimensional. This yields a curve within the connected subvariety $S_{n-k}(\Phi'|_c) \subset Gr_{n-k}(E'|_c)$ of vector subspaces of the fiber $E'|_c$ which are invariant by $\Phi'|_c$. As argued in the proof of \cite[Theorem 4.16]{HauselHitchin}, this translates into a curve in the moduli space of Higgs bundles, defined by $\mathcal H_{V_\lambda}(E', \Phi') \simeq \lambda \mathcal H_V(E',\Phi')$. This curve connects $\mathcal H_{V_0}(E',\Phi') = (E(-c),\Phi)$ with a different fixed point, given by $\mathcal H_{V_\infty}(E',\Phi')$, where $V_\infty = \gen{v_{k-1}, v_k, v_{k+1}, \dots,v_{k+l-1},v_{k+l+1}, v_{n-1}}$. This fixed point is also stable \cite[Lemma 4.18]{HauselHitchin}.
	
	Moreover, since $l$ is odd, then $l+1$ is even and we have that $V_{\lambda} = V_{-\lambda}$, hence $\lambda \mathcal H_V(E',\Phi') \simeq -\lambda \mathcal H_V(E',\Phi')$ so that this curve is fixed by the action of $-1$ and hence it lies in $W^{2+}_{(E(-c),\Phi)}$. This shows that $(E(-c), \Phi)$ is even wobbly thus $(E, \Phi)$ as well. \qed

There are still more examples of even wobbly stable Higgs bundles of type $(1,1,\dots,1)$, which are covered by the following proposition.

\begin{proposition}\label{wobbly2}
    Let $\mathcal E = (E, \Phi)$ be a smooth fixed point of rank $n \ge 4$. Suppose that $b := b_{n-2} \circ \dots \circ b_2$ has a multiple zero at $c \in C$. Then $\mathcal E$ is even wobbly.
\end{proposition}

    \noindent {\it {\it Proof.}} The proof is identical to that of Proposition \ref{wobbly1}, with the only difference being the construction of the invariant subspace $V \in S_{n-k}(\Phi'|_c)$. 
    
    Let $1 < i_1 < i_2 < \dots < i_m < n-1$ be all the indices other than $1$ and $n-1$ such that $b_{i_k}(c) = 0$. We can assume that all of them are of the same parity since otherwise we apply the previous proposition. We start again with a Hecke transform at $(L_0 \oplus L_1 \dots \oplus L_{i_m-1})|_c$. This yields $(E', \Phi')$ with the same form as in the proof of Proposition \ref{wobbly1}.

Now we will construct a basis $\{v_0,v_1, \dots, v_{n-1}\}$ of $E'|_c$ as follows. We start with any nonzero $v_0 \in L_0'|_c$, then apply $\Phi'|_c$ until a zero vector is obtained. This will produce a \textit{string}: $v_1 = \Phi'|_c(v_0)$, $v_2 = \Phi'|_c(v_1) \dots$ all the way up to $v_{i_1-1}$, if $b_1(c) \neq 0$, or just $v_0$ if $b_1(c) = 0$. The process then iterates: we pick $v_{i_1}$ (or $v_1$ if $b_1(c) = 0$) nonzero in the corresponding summand and iterate, repeating until $v_{n-1}$. By construction, the basis is partitioned in strings inside of which each element maps to the next one and the last one maps to zero. Also, each $v_i$ is in the corresponding summand $L'_i|_c$. Let $V_0 = \gen{v_{i_m}, \dots, v_{n-1}}$.

Now recall that $i_m < n-1$ and $i_1 \ge 2$. Note that if $m = 1$, it is possible that $i_1 = i_m$, meaning a multiple zero of $b_{i_1}$ at $c$. In any case, we define the following $(n-i_m)$-dimensional subspace:
$$V = \gen{v_{i_m} + v_{i_1-2}, v_{i_m+1}+v_{i_1-1}, v_{i_m+2}, \dots, v_{i_m+i_1-1}, v_{i_m+i_1}, \dots, v_{n-1}}.$$

It is $\Phi'|_c$-invariant: each generator is taken to the next, and the last one to zero. Notice that it is important that $i_m < n-1$ since otherwise we get a $1$-dimensional space that does not work, for $\gen{v_{i_m} + v_{i_1-2}} \mapsto \gen{v_{i_1-1}}$. Also notice that if $b_{n-1}(c) = 0$ then the second generator might map to zero instead of $v_{i_m+2}$, but this is not a problem for the invariance. We compute $\lambda V$:
$$V = \gen{v_{i_m} + \lambda^{i_m-i_1+2}v_{i_1-2}, v_{i_m+1}+\lambda^{i_m-i_1+2}v_{i_1-1}, v_{i_m+2}, \dots, v_{i_m+i_1-1}, v_{i_m+i_1}, \dots, v_{n-1}}.$$

Once again this gives a curve in $S_{n-k}(\Phi'|_c)$ connecting $V_0$ with a different $V_\infty$ such that the Hecke transform is a new fixed point and, since $i_m - i_1 + 2$ is even, it follows that $\lambda V = -\lambda V$, as desired. \qed

We now prove that these constitute all even wobbly cases.

\begin{proposition}
	Let $\mathcal E_\delta = (E, \Phi)$ be a smooth fixed point of type $(1,1,\dots,1)$. Suppose that for every $c \in C$ we have that $b_i(c) = b_j(c) = 0$ implies $i \equiv j \mod 2$, and $b_{n-2} \circ \dots \circ b_2$ has at most a single zero at $c$. Then, $\mathcal E$ is even very stable.
\end{proposition}

\noindent {\it Proof.} We argue by induction on $\deg(\delta_{n-1} + \dots + \delta_1)$, as we already know from Section \ref{examples} that the degree $0$ case, which is the canonical uniformising Higgs bundle, is even very stable. Suppose that $(E',\Phi') \in W^+_{\mathcal E} \cap \mathcal C$ is a nilpotent element in the upward flow of $\mathcal E$, with the full filtration
$$0 = W_0 \subset W_1 \subset \dots \subset W_{n-1} \subset W_n = E'$$ 

\noindent given by Proposition \ref{upwardflow}. Take one of the points $c \in C$ with $\mu_\delta(c) \neq 0$, and suppose first that $b_{n-1}(c) = 0$. This means that $V' := (W_{n-1})|_c$ is $\Phi'|_c$-invariant, and, as explained in \cite[Section 4]{HauselHitchin}, the Hecke transformation $(E'_1,\Phi'_1) = \mathcal H_{V'}(E',\Phi')$ is stable, nilpotent and lies in the upward flow of $(E_1,\Phi_1) := \mathcal H_V(E,\Phi)$, where $V=(L_0 \oplus \dots \oplus L_{n-2})|_c \subseteq E|_c$. We have that $(E_1,\Phi_1)$ is another type $(1,1,\dots,1)$ point given by $\delta^1 = (\delta_0^1,\dots,\delta_{n-1}^1)$ with the only difference being that $\delta_{n-1}^1 = \delta_{n-1} - c$. Thus, by the induction hypothesis it follows that $(E_1,\Phi_1)$ is even very stable, which results in the following two options:

\begin{itemize}
	\item If $(E_1',\Phi_1') \not\simeq (E_1,\Phi_1)$, since the latter is even very stable, we have that $(E_1',\Phi_1')$ lies in an odd weight space of $T_{(E_1,\Phi_1)} \mathcal M$ and, by using the invariant subspaces $V_1' \subseteq E_1'|_c$ and $V_1 \subseteq E_1|_c$ such that $\mathcal H_{V_1'} (E_1',\Phi_1') = (E',\Phi')$ and $\mathcal H_{V_1} (E_1,\Phi_1) = (E,\Phi)$ (up to twisting by a fixed line bundle), it follows that $(E',\Phi')$ also has odd weight in the upward flow for $(E,\Phi)$ and hence $(E',\Phi') \notin W_{\mathcal E_{\delta}}^{2+}$.
	\item If $(E_1',\Phi_1') \simeq (E_1,\Phi_1)$, then the one dimensional subspace $V_1'\subseteq E_1'|_c$ such that $\mathcal H_{V_1'} (E_1',\Phi_1') = (E',\Phi')$ is one of the $\Phi_1|_c$-invariant subspaces in the pair $(E_1,\Phi_1)$, which we know well since it comes from a Hecke transformation of the starting $(E,\Phi)$. Indeed, we have:
	$$E_1 = L_0 \oplus \dots \oplus L_{n-2} \oplus L_{n-1}(-c),$$
	
	$$\Phi_1 = \begin{pmatrix}
		0 & 0 & \dots & 0 &0& 0 \\
		b_1 & 0 & \dots & 0  &0& 0 \\
		0 & b_2 & \dots & 0 &0& 0 \\
		\vdots & \vdots & \ddots & \vdots & \vdots & \vdots \\
		0 & 0 & \dots & b_{n-2}  & 0& 0\\
		0 & 0 & \dots & 0  & \frac{b_{n-1}}{s_c} & 0
	\end{pmatrix}.$$

	Let $V = \{v_0,\dots,v_{n-1}\}$ be a basis of $E_1|_c$. From the starting hypotheses about $(E,\Phi)$, the only basis vectors that could vanish via $\Phi_1|_c$ are $\{v_0, v_{j-1}, v_{n-2}, v_{n-1}\}$ where $j$ is the only index with $2 \le j \le n-2$ such that $b_j(c) = 0$, if it exists. Necessarily, $j \equiv n-1 \mod 2$. Also, if $v_0$ vanishes then $b_1(c) = 0$ and hence $n \equiv 0 \mod 2$. Finally, $v_{n-2}$ vanishes if and only if $b_{n-1}$ had a multiple zero at $c$.
	
	Since $\Phi_1$ is nilpotent, the desired invariant $1$-dimensional subspace must be of the form $V_1' = \gen{v}$ where $\Phi_1|_c(v) = 0$, that is, $v=\alpha v_0+\beta v_{j-1}+\gamma v_{n-2}+ \delta v_{n-1}$. Moreover, 
	$$\lambda \cdot V_1' = \gen{\alpha \lambda^{n-1} v_0+\beta \lambda^{n-j} v_{j-1}+\gamma \lambda v_{n-2}+ \delta v_{n-1}}$$
	 
	must verify that $\lim\limits_{\lambda \to 0}\lambda V_1' = V_1 = \gen{v_{n-1}}$, which implies that $\delta \neq 0$. Hence, because of the parity of the exponents of $\lambda$ appearing in the expression for $\lambda \cdot V_1'$, if $(E',\Phi') = \mathcal H_{V_1'}(E_1',\Phi_1')$ were in an even weight space, that is, if $\lambda \cdot V_1' = -\lambda \cdot V_1'$, it would be necessary that $\alpha=\beta=\gamma=0$, thus $V_1' = V_1$, meaning $(E',\Phi') \simeq (E,\Phi)$.
\end{itemize}

This concludes the analysis for the $b_{n-1}(c) = 0$ case. Now, if $b_{n-1}(c) \neq 0$, we may also assume that $b_1(c) \neq 0$. This is because the involution $(E,\Phi) \mapsto (E^*,\Phi^t)$ of $\mathcal M$ naturally bijects fixed points and upward flows, sending a point with $b_1(c) = 0$ to a point with $b_{n-1}(c) = 0$. Hence, the remaining case is when the only $b_j(c) = 0$ happens at a single $j$ with $2 \le j \le n-2$ and multiplicity one. This case is treated exactly as in the proof of \cite[Theorem 4.16]{HauselHitchin}, which we now recall.

We have the $\Phi'|_c$-invariant $j$-dimensional subspace $V' := (W_{j-1})|_c$ giving a nilpotent $(E_1',\Phi_1') = \mathcal H_{V'}(E',\Phi')$ in the upward flow of $(E_1,\Phi_1) = (E_1,\Phi_1) := \mathcal H_V(E,\Phi)$, where $V=(L_0 \oplus \dots \oplus L_{j-1})|_c \subseteq E|_c$. Exactly as before, if $(E_1',\Phi_1') \not\simeq (E_1,\Phi_1)$, the induction hypothesis gives $(E',\Phi') \notin W_{\mathcal E_\delta}^{2+}$. Otherwise, notice that $\Phi_1'|_c$ is now a regular nilpotent, so the only $\Phi_1'|_c$-invariant $(n-j)$-dimensional subspace is $V_1' = V_1 = (L_{j} \oplus \dots \oplus L_{n-1})|_c$. Hence (up to tensoring everything by a fixed line bundle) we have $(E',\Phi') = \mathcal H_{V_1}(E_1,\Phi_1) = (E,\Phi)$. \qed

\subsection{Even minuscule weights in $\GL(n,\mathbb C)$}
We will now see how the conditions found before arise naturally in the context of the root system of the structure group $\GL(n,\mathbb C)$ of $E$. We have $\omega_1, \dots, \omega_{n-1}, \omega_n\in\Lambda^+$  the fundamental weights.  Recall that the {\em height} of a root is the number of simple roots in its decomposition. In terms of the fundamental weights	the positive roots of height $k\in \{1,\dots,n-1\}$  can be indexed by $p \in \{1,\dots,n-k\}$ and given by: $$\alpha_{k,p}:=-\omega_{p-1} + \omega_{p}+\omega_{p+k-1} - \omega_{p+k},$$ where $\omega_0$  is understood as $0$. In particular, the simple roots are the positive roots of height $1$ of the form $\alpha_{1,p}=-\omega_{p-1} + 2\omega_{p} - \omega_{p+1}$ for $p=1,\dots,n-1$. The highest root  $\alpha_{n-1,1}=\omega_1+\omega_{n-1}-\omega_n$ has height $n-1$.


\begin{definition}
	For dominant weights $\lambda, \mu$, we define the {\em even partial ordering} as $$\lambda \ge_2 \mu \iff \lambda - \mu \in \Phi^{2+},$$ where $$\Phi^{2+} = \left\{\sum_{\substack{k \text{ even,}\\ 1 \le p \le n-k}}c_{k,p}\alpha_{k,p} : c_{k,p} \in \mathbb Z_{\ge 0}\right\}$$ is the set of positive linear combinations of even positive roots. Minimal elements for this ordering are called {\em even minuscule}.
\end{definition}

We will say for a weight $\lambda$ that its \textit{$i$-th coordinate/entry}, $\lambda_i$, is the coefficient of $\omega_i$ when $\lambda$ is written in the basis $\{\omega_1, \dots, \omega_n\}$. The \textit{position} of the \textit{$i$-th} coordinate will be just $i$. Now we will characterize the even minuscule weights, seeing that the conditions exactly match those for even very stable Higgs bundles. The characterization will be carried out solely via combinatorial arguments.

\begin{proposition}
	Let $\lambda$ be a dominant weight such that at least one of these holds:
	
	\begin{enumerate}
		\item The weight $\lambda$ has nonzero coordinates at two positions $1 \le i < j \le n-1$ with $i \not\equiv j \mod 2$.
		\item The weight $\lambda$ verifies $\lambda_2 + \dots + \lambda_{n-2} \ge 2$.
	\end{enumerate}
	
	\noindent Then $\lambda$ is not even minuscule.
\end{proposition}

\noindent {\it Proof.} For the first situation, if $\lambda$ has nonzero coordinates at $i$ and $j$, where $j - i := k$ is an odd positive number, then consider $\mu := \lambda - \alpha_{i,k+1}$. Because of the nonzero coordinates at $i$ and $j$ it follows that $\mu$ is dominant. By construction, $\mu < \lambda$.

If $\lambda_2 + \dots + \lambda_{n-2} \ge 2$, choose indices $i,j$ such that $2 \le i \le j \le n-2$ with either $i \neq j$ and $\lambda_i, \lambda_j \ge 1$, or $i = j$ and $\lambda_i \ge 2$. We can assume the positions are of the same parity, since otherwise we apply the previous. Hence write $k := j-i$ a nonnegative even number. Defining $\mu := \lambda - (\alpha_{k+2,i}+\alpha_{k+2,i-1})$ works as before, since $\alpha_{k+2,i}+\alpha_{k+2,i-1} = -\omega_{i-1}+\omega_i+\omega_{j+1}-\omega_{j+2} - \omega_{i-2} + \omega_{i-1} + \omega_{j} - \omega_{j+1} = -\omega_{i-2} + \omega_i + \omega_j -\omega_{j+2}$. \qed

\begin{remark}
	Notice in the last part of the previous proof that $2 \le i,j \le n-2$ is indeed required: for example, if $i = 1$ we cannot take $p = i-1 = 0$. Similarly we need for $\alpha_{k+2,i}$ to make sense that $i \le n-(k-2) = n-k+2$ which is equivalent to $j \le n-2$.
	
	As an example for the proof, take $n=9$ and $ \lambda = (0,1,0,0,0,1,0,0,0)$. Following the proof, we can consider $\alpha_{6,2}+\alpha_{6,1} = (-1,1,0,0,0,0,1,-1,0) + (1,0,0,0,0,1,-1,0,0) = (0,1,0,0,0,1,0,-1,0)$, so that $\mu = (0,0,0,0,0,0,0,1,0)$ is lower than $\lambda$.
\end{remark}

 For the reciprocal, the combinatorial arguments will be easier if we stop considering $\omega_0$ as zero and rather see it as an extra linearly independent vector. That is, we view the weight space $W:=\gen{\omega_1,\dots,\omega_{n-1}, \omega_n}$ as a subspace of a new vector space $\widetilde{W}:=\gen{\omega_0,\omega_1,\dots,\omega_{n-1}, \omega_n}$ which is a dimension higher. We have a projection $\pi : \widetilde{W} \to W$. We define $\widetilde{\alpha}_{k,p} := -\omega_{p-1} + \omega_{p}+\omega_{p+k-1} - \omega_{p+k}$ as before but in $\widetilde{W}$ (i.e. we do not consider $\omega_0=0$ in that expression anymore), so that $\pi(\widetilde{\alpha}_{k,p}) = {\alpha}_{k,p}$. We also lift the even positive lattice as:

$$\widetilde{\Phi}^{2+} = \left\{\sum_{\substack{k \text{ even,}\\ 1 \le p \le n-k}}c_{k,p}\widetilde{\alpha}_{k,p} : c_{k,p} \in \mathbb Z_{\ge 0}\right\},$$

\noindent so that $\pi (\widetilde{\Phi}^{2+}) = \Phi^{2+}$.

For example, in this new vector space the simple positive roots have the following lifts:  $\widetilde{\alpha}_1 = (-1,2,-1,0,0,\dots)$, $\widetilde{\alpha}_2 = (0,-1,2,-1,0,0,\dots)$ etcetera. As can be seen, this avoids the situation of the negative coordinate being truncated at the beginning. In this setting we have the following facts about $\widetilde{\Phi}^{2+}$, which will be the only ones we will need for our proof:

\begin{lemma}\label{combinatorics}
	Take a nonzero $\widetilde{x} = (x_0,x_1,\dots,x_{n-1},x_n) \in \widetilde{\Phi}^{2+}$. We have:
	
	\begin{enumerate}
		\item The values $x_0$ and $x_n$ are not positive.
		\item The values $\sum\limits_{j=0}^nx_j$, $\sum\limits_{j \text{ even}} x_j$ and $\sum\limits_{j \text{ odd}} x_j$ are all zero.
		\item The value $\sum\limits_{j=2}^{n-2}x_j$ is not negative.
		\item The first and last nonzero coordinates of $\widetilde{x}$ are negative.
	\end{enumerate}
\end{lemma}

\noindent {\it Proof.} Immediate, by induction on the number of positive even roots into which $\widetilde{x}$ decomposes. First, it is clear that any $\widetilde{\alpha}_{k,p}$ for even $k$ has all those properties. Second, it is easy to check that the sum of any two vectors with those properties keeps satisfying them. Hence the result follows. \qed

\begin{proposition}
	Let $\lambda$ be a dominant weight that is not even minuscule. Then at least one of these hold:
	
	\begin{enumerate}
		\item The weight $\lambda$ has nonzero coordinates at two positions $1 \le i < j \le n-1$ with $i \not\equiv j \mod 2$.
		\item The weight $\lambda$ verifies $\lambda_2 + \dots + \lambda_{n-2} \ge 2$.
	\end{enumerate}
\end{proposition}

\noindent {\it Proof.} Take $\mu < \lambda$. Denote $x := \lambda - \mu \in \Lambda^{2+}$. Since $\lambda = x+\mu$ and the coordinates of $\mu$ at positions between $1$ and $n-1$ are non negative, it suffices to check that either $x$ has positive coordinates at positions between $1$ and $n-1$ of different parity, or that the sum of the positive coordinates in positions $\{2,\dots,n-2\}$ of $x$, which we shall denote $S$ from now on, verifies $S \ge 2$. 

In order to do this we can work with a lift $\widetilde{x} \in \widetilde{\Lambda}^{2+}$ of $x$ by lifting each of the positive even roots in some decomposition. We will prove that $\widetilde{x}$ has at least one of the two desired properties, and, since $x_0,x_n \le 0$ by Lemma \ref{combinatorics}, then $x = \pi(\widetilde{x})$ will also have them and the proof will be complete.

Assume that $\tilde{x}$ has every positive coordinate between $1$ and $n-1$ at even positions. We will show that $S \neq 0$ and $S \neq 1$, hence $S \ge 2$.

If $S = 0$ this means there are no positive entries in positions $\{2,\dots,n-2\}$. By Lemma \ref{combinatorics}, part $3$ there are no negative entries either. So we have $\tilde{x} := (x_0,x_1,0,\dots,0,x_{n-1},x_n)$. Now if $n$ is even then by Lemma \ref{combinatorics}, part $2$ we have $x_0+x_n = 0$ hence $x_0 = x_n = 0$, the latter by Lemma \ref{combinatorics}, part $1$. But, since $x_1 = -x_{n-1}$ this contradicts Lemma \ref{combinatorics}, part $4$ reaching a contradiction. If $n$ is odd then $x_0+x_{n-1} = 0$ and $x_1 + x_n = 0$. By assumption that positive entries are in even positions, we have $x_1 = x_n = 0$ and we reach the same contradiction.

If $S = 1$ there will be exactly one positive entry in positions $\{2,\dots,n-2\}$ of value $1$, say the $i$-th one. By assumption $i$ is even. By Lemma \ref{combinatorics}, part $3$, from positions $2$ to $n-2$ there could either be no negative entries or one negative entry of value $-1$ at position $j$. Now a case by case analysis follows.

\begin{itemize}
	\item If $n$ is odd, the case $j$ odd is excluded because we would get, summing odd entries, that $x_1 + x_n - 1 = 0$ and we have $x_1, x_n \le 0$, the first by assumption and the second by Lemma \ref{combinatorics}, part $1$. So we have $x_1=x_n= 0$ and either $j$ is even or there is no $j$ at all. In either case we find, summing the even entries, that $x_0+x_{n-1} \in \{0,-1\}$, from which it is impossible that both are negative, contradicting Lemma \ref{combinatorics}, part $4$.
	\item If $n$ is even, the case $j$ odd is again excluded since, summing odd entries, we get $x_1+x_{n-1} - 1 = 0$, but by assumption both $x_1$ and $x_{n-1}$ are not positive. Hence $x_1 = x_{n-1} = 0$ and $j$ is even or there is no $j$ at all. In this case we have to distinguish: if there is no $j$ then $x_0+x_n = -1$ so both cannot be negative at the same time, and since the only other nonzero entry is the $1$ at the $i$-th position, Lemma \ref{combinatorics}, part $4$ is contradicted. If there is a $j$ then $x_0 + x_n = 0$, so that $x_0=x_n=0$ and once again we have a contradiction with Lemma \ref{combinatorics}, part $4$.
\end{itemize}

The only remaining task to complete the proof is to work out the case where all positive entries of $\tilde{x}$ are at odd positions instead of even ones. There is a symmetry in $\tilde{\Phi}^{2+}$ via $\omega_{k} \mapsto \omega_{n-k}$ that will change odd positions and even positions if $n$ is odd. So the only remaining cases are for even $n$, which can be approached exactly the same as before:

\begin{itemize}
	\item If $S = 0$ then the previous argument for even $n$ did not use the positions of positive entries being even so it still works when they are odd.
	\item If $S = 1$ we now have that $i$ is odd. If $j$ were even, summing the even entries we get $x_0+x_n-1 = 0$ which is not possible ($x_0,x_n \le 0$) so that $j$ is either odd or does not exist. If it is odd, we sum the odd entries to get $x_1+x_{n-1} = 0$, which yields a contradiction with Lemma \ref{combinatorics}, part $4$. If it does not exist, we get $x_1+x_{n-1} = -1$, from which they cannot both be negative and we reach the same contradiction.
\end{itemize}
\qed

Thus, we can rephrase the condition for even wobbliness as follows:

\begin{theorem}\label{evenmin}
	Let ${\mu}:C\to \Lambda^+$. Then $\calE_{\mu}\in \M^{s\T}$  is even very stable if and only if for every $c \in C$, the dominant weight ${\mu}(c)$ is even minuscule.
\end{theorem}

As an example, the fixed point in rank $n=4$ associated to $\mu(c) = \omega_1+\omega_3$ and $\mu(d) \neq 0$ for $d \neq c$ is even very stable, while being wobbly in the usual sense. In rank $n=2$ every fixed point is even very stable, which can be deduced more easily by the lack of even weights in downward flows. This is because the only weights appearing in the upward flows in rank $n=2$ are those appearing in the Hitchin base, which are $1$ and $2$. We then know that a negative weight $k$ is paired to a nonnegative one, $1-k$, via the symplectic form $\omega$, hence the only possible negative weight is $k=-1$.
 
As a final remark, it is possible to work with higher order automorphisms defined by $(E,\Phi) \mapsto (E,\zeta_r\Phi)$, where $\zeta_r$ is a primitive $r$-th root of unity, giving subspaces $\mathcal M^r \subset \mathcal M$ and the corresponding notions of $r$-very stable bundles. Most of the results presented here naturally generalize to that situation (by using $r$-minuscule weights), however in this case there is no longer an associated real group and the subspaces considered are not Lagrangian.

\section{Hitchin map on upward flows}
\label{Hitchin}
One important motivation to consider very stable upward flows in \cite{HauselHitchin} was the observation that  the Hitchin map $$h_\calE:=h|_{W^+_\calE}:W^+_\calE\to \calA$$ restricted to them is proper. Furthermore $W^+_\calE\cong T_\calE^+\calM$ as $\T$-varieties. Thus in the very stable case $h_\calE$ is a proper, even finite flat \cite[Lemma 4.6]{HauselHitchin}, $\T$-equivariant morphism between semi-projective affine spaces of the same dimension. As such it is suspectible for explicit description. 

\subsection{Equivariant cohomology of homogeneous spaces}

We will describe the Hitchin map explicitely on some very stable upward flows in terms of equivariant cohomology. First we recall some of the basic properties of equivariant cohomology see e.g. \cite{atiyah-bott} for more details. 

Let $\G$ be a connected complex affine (or compact Lie) group. 
Consider the classifying $\G$-bundle $E\G\to B\G$, where $E\G$ is contractible. $B\G$ is called the classifying space, and its cohomology ring can be computed as follows: $$H^*_\G:=H^*(B\G;\C)\cong H^*(B{\rm T};\C)^{\W_\G}\cong \C[\t]^{\W_\G},$$ where $\t=\op{Lie}({\rm T})$ is the Lie algebra of a maximal torus ${\rm T}\subset \G$ and $\W_\G$ is the Weyl group of $\G$. 

Let now $\G$ act on a variety (or manifold) $X$. We can form the Borel, or homotopy, quotient $X_\G:=(E\G\times X)/\G$, by the diagonal action of $\G$, which is an $X$-bundle over $B\G$. Its cohomology ring is what is taken to be its equivariant cohomology $$H^*_\G(X;\C):=H^*(X_\G;\C).$$ As $X_\G$ is an $X$-bundle over $B\G$ we get a ring map $H^*_\G\to H^*_\G(X;\C)$, making $H^*_\G(X;\C)$ an $H^*_\G$-algebra. 

Let now ${\rm T}<{\rm H}<\G$ be a connected closed subgroup containing the maximal torus ${\rm T}$ of $\G$. Because $E\G$ is contractible, we have $E\G/{\rm{H}}\sim E{\rm{H}}/{\rm{H}}$ so we can compute the $\G$-equivariant cohomology of the homogeneous space $\G/{\rm{H}}$ as follows: $$H^*_\G(\G/{\rm{H}};\C)\cong H^*((E\G\times \G/{\rm{H}})/\G;\C)\cong H^*(E\G/{\rm{H}};\C)\cong H^*_{{\rm{H}}},$$  with the ring map $H^*_\G\to H^*_{{\rm{H}}}$ from the natural map $B{\rm{H}}\to B\G$ induced from the embedding ${\rm H} \subset \G$. This way we have a simple way to compute equivariant cohomology of a homogeneous space explicitly as in the following diagram of graded algebras: 
\begin{align}\label{eqcoh} \begin{array}{ccc}
    H^*_\G(\G/{\rm{H}};\C) & \cong  & \C[\t]^{\W_{\rm{H}}} \\
    \uparrow & &  \uparrow\\ 
    H^*_\G & \cong & \C[\t]^{\W_{\G}}.
\end{array}
\end{align} Note that the structure maps above are finite free, meaning that they define a finite free module. This is because the equivariant cohomology of equal rank homogeneous spaces are equivariantly formal \cite[(1.2)]{goresky-kottwitz-macpherson}, for example because they have no odd cohomology \cite[Theorem VII, p. 467]{greub-etal}. 

We will see below that on certain very stable upward flows the Hitchin map can be modelled by the spectrum of equivariant cohomology of appropriate homogenous spaces. For this reason we also take the spectrum of \eqref{eqcoh} and record the corresponding diagram:
\begin{align}\begin{array}{ccc}
   \Spec(H^*_\G(\G/{\rm{H}};\C)) & \cong  & \t/\!\!/{\W_{\rm{H}}} \\
    \downarrow & &  \downarrow\\ 
    \Spec(H^*_\G) & \cong & \t /\!\!/ {\W_{\G}}.
\end{array} 
\end{align} We note that half the grading on cohomology will induce the $\T$-action on $\t/\!\!/{\W_{\rm{H}}}$, which agrees with the $\T$-action induced from  weight one action on $\t$. The down arrows in the diagram then become $\T$-equivariant finite flat (ultimately because of equivariant formality), in particular proper, morphisms. 

Because we will also consider how certain involutions act on the equivariant cohomology of homogeneous spaces, here we record the following lemma. To formulate it recall that for a unital commutative ring $R$ with $2$ invertible and with involution $\theta:R\to R$, the {\em coinvariant ring} is defined as $$R_\theta:=R/(r-\theta(r))_{r\in R}\cong R/(r \in R : \theta(r)=-r).$$ While if $A$ is a commutative $R$-algebra and $\theta:A\to A$ also acts on $A$, compatibly with the action on $R$, then $A_\theta$ is naturally an $R_\theta$-algebra, called the {\em coinvariant algebra}. Their relevance is in forming the fixed point scheme of the affine $R$-scheme $\Spec(A)$ under the involution $\Spec(\theta):\Spec(A)\to \Spec(A)$, which we also denote by $\theta$. We get that the fixed point scheme \begin{align}\label{fixed}\Spec(A)^\theta\cong \Spec(A_\theta)\end{align} is an affine $R_\theta$-scheme. Using arguments from \cite{PanyushevYakimova} we have the following 

\begin{lemma}\label{taustarcoinvariant} Let $\G$  be a connected complex reductive group (or a  connected compact Lie group), and $\tau:\G\to \G$ a complex algebraic (smooth) involution, with ${\rm T}\subset \G$ a $\tau$-stable torus. Let $\G_0^\tau$ denote the identity component of the fixed point group $\G^\tau$. Assume that \begin{align}
    H^*_\G\twoheadrightarrow H^*_{\G^\tau_0} \label{surject} \end{align} is surjecive.  Let $\theta:=\tau^*:H^*_\G\to H^*_\G$ denote the induced action on the cohomology of $B\G$.  Then we have that the coinvariant algebra $$(H^*_\G)_\theta\cong H^*_{\G_0^\tau}.$$ Moreover if ${\rm T}\subset {\rm H}\subset\G$ is a $\tau$-invariant closed connected subgroup, such that $H^*_{\rm H}\twoheadrightarrow H^*_{{\rm H}^\tau_0}$ is surjective, then we have the commutative diagram: $$\begin{array}{ccccccc} H^*_\G(\G/{\rm H};\C)_\theta&\cong&(H^*_{\rm{H}})_\theta &\cong & H^*_{\rm{H}^\tau_0} &\cong&H^*_{\G_0^\tau}(\G_0^\tau/{\rm H}_0^\tau;\C) \\ \uparrow &&\uparrow&& \uparrow & &\uparrow \\ (H^*_{\G})_\theta  &\cong& (H^*_{\G})_\theta  &\cong& H^*_{\G_0^\tau}&\cong& H^*_{\G_0^\tau}\end{array}$$
    \end{lemma} 
 \noindent   {\it Proof.}  We know that $H^*_\G\cong \C[\t]^\W$, the algebra of invariant polynomials on $\t$ by the Weyl group $\W$ of $\G$. By \cite[Lemma 6.1]{Springer} we can choose  algebra generators $p_1,\dots,p_r$ in $\C[\t]^\W$ such $\tau(p_i)=\epsilon_i p_i$, where $\epsilon_i=\pm 1$. By  \cite[Lemma 6.5]{Springer} the generators with $\epsilon_i=1$ give the generators of $\C[\t^\tau]^{\W^\tau},$ where $W^\tau$ is the subgroup of $W$ fixed by $\tau$, itself a reflection group on $\t^\tau$. Thus the number of $p_i$'s with $\epsilon_i=1$ is exactly the dimension of $\t^\tau$, which in turn agrees with the rank of ${\rm \G_0^\tau}$. 
 
 Finally, all $p_i$'s with $\epsilon_i=-1$ restrict trivially to $\t^\tau$. Thus by the assumption \eqref{surject}  the restriction of those with $\epsilon_i=1$ generate $\C[\t^\tau]^{\W_{\G_0^\tau}}\cong H^*_{\G_0^\tau}$, and as there are $\rank(\G_0^\tau)$ of them, they should restrict algebraically independent. Thus the kernel of the surjection is generated by the anti-invariant generators, showing the claim $(H^*_\G)_\theta\cong H^*_{\G_0^\tau}$.

 The second statement follows from the first and \eqref{eqcoh}. \qed

    \begin{remark}\label{coinvariantsymmetric}
    The assumption of surjection \eqref{surject} is quite restrictive. For simple $\G$ it only happens \cite[after (4.1)]{PanyushevYakimova} for the symmetric pairs $$(\G,\G^\tau_0)\cong (\SL_{2n+1},\SO_{2n+1}), (\SL_{2n},\Sp_n), (\SO_{2n},\SO_{2n-1}), \mbox{ and } (\rm{E}_6,\rm{F}_4).$$ All these examples will appear in \S\ref{thetamodell} below.  
    \end{remark}

\subsection{Explicit Hitchin map on  very stable upward flows}
\label{modell}

First we recall from \cite{Hausel-hmsec} how to describe $h_k:=h_{\calE_k}$ for the very stable upward flows $W_k^+:=W^+_{\calE_k}$ explicitly, where $\calE_k:=\calE_{\mu_c^k}$, $c\in C$ fixed and $\mu^k_c:C\to \Lambda^+$ is defined by $\mu^k_c(c)=\omega_k$ and $0$ otherwise. 

Let us define the evaluation map at $c\in C$  $${\rm ev}_c:\calA\to \C^n$$ sending the characteristic polynomial $(a_{1},\dots,a_n)\in \calA$ to $$(a_1(c),\dots,a_n(c))\in K^1_c\times \dots \times K^n_c\cong \C^n$$ after identifying the fiber $K\cong \C$. We can further identify $$\C^n\cong \gl_n/\!\!/\GL_n\cong \t_n/\!\!/{S_n}\cong \Spec(H^*(B{\GL}_n,\C)),$$ where $\gl_n={\op{Lie}}(\GL_n)$ and $\t_n={\op{Lie}}({\rm T}_n)$ the Lie algebra of the maximal torus ${\rm T}_n\subset \GL_n$, and the symmetric group $S_n$ is the Weyl group of $\GL_n$. Then with the notation $H^*_{\GL_n}:=H^*(B{\GL}_n;\C)$ for the cohomology ring of the classifying space $B{\GL}_n$ we have the following pull back diagram from \cite{Hausel-hmsec}

$$\begin{array}{ccc} W^+_k & \twoheadrightarrow & \Spec(H^{2*}_{\GL_n}(\rm{Gr}_k(\C^n),\C))\\\!\!\!\!\!\! h_k \downarrow & \lrcorner & \downarrow \\ \calA &\stackrel{{\rm ev}_c}{\twoheadrightarrow} & \Spec(H^{2*}_{\GL_n}) \end{array},$$ where $\Gr_k(\C^n)$ is the Grassmannian of $k$-dimensional subspaces in $\C^n$ with the usual action of $\GL_n$. Additionally, all maps are $\T$-equivariant with respect to the usual $\T$-action on $W_k^+\subset \M$ and $\calA$, and the one induced by the grading on $H^{2*}_{\GL_n}(\Gr_k(\C^n);\C)$ and $H^{2*}_{\GL_n}$. 
In other words the Hitchin map $h_k$ can be modelled by the equivariant cohomology of the Grassmannian ${\rm Gr}_k(\C^n)$. 

\subsection{Explicit Hitchin map on even very stable upward flows}
\label{thetamodell}

We can then ask what models the even Hitchin map $h^\theta_{k}: \W_k^{2+}\to \calA^\theta$. We can note that $\theta$ acts on $\M$ as $-1\in \T$ in the natural $\T$-action, thus we can induce an action of $\theta$ on $H^{2*}_{\GL_n}(\Gr_k(\C^n);\C)$ and $H^{2*}_{\GL_n}$ as $-1\in \T$. Therefore we have 
\begin{align} \label{thetamodel} \begin{array}{ccc} W^{2+}_k & \twoheadrightarrow & \Spec(H^{2*}_{\GL_n}(\rm{Gr}_k(\C^n),\C))^\theta \\ \!\!\!\!\!\! h^\theta_k \downarrow & \lrcorner & \downarrow \\ \calA^\theta &\stackrel{{\rm ev}_c}{\twoheadrightarrow} & \Spec(H^{2*}_{\GL_n})^\theta \end{array},\end{align}
and, as a result, the even Hitchin map $h^\theta_k$ can be modelled on the $\theta$-fixed point scheme $$\Spec(H^{2*}_{\GL_n}(\rm{Gr}_k(\C^n),\C))^\theta.$$

For simplicity we will start with $\GL_{2n}$ and an even number $0<2k<2n$. To understand the fixed point scheme $\Spec(H^{2*}_{\GL_n}(\rm{Gr}_{2k}(\C^{2n}),\C))^\theta$ we recall a presentation of the equivariant cohomology of the Grassmannian. In practice it can be done by following through the restriction of invariant polynomials 
in \eqref{eqcoh}.

Let $e_1,\dots,e_{2k}, f_1,\dots,f_{2n-2k}, c_1,\dots,c_{2n}$ be variables of degree given by their index. Then we have the following presentation of the graded ring $H^{2*}_{\GL_n}(\Gr_{2k}(\C^{2n}),\C)\cong$ $$ \frac{\C[e_1,\dots,e_{2k}, f_1,\dots,f_{2n-2k}, c_1,\dots,c_{2n}]}{\left( (t^{2k}+e_1t^{2k-1}+\dots +e_{2k})(t^{2n-2k}+f_1t^{2n-2k-1}+\dots +f_{2n-2k})-(t^{2n}+c_1t^{2n-1}+\dots+c_{2n})\right)},$$ where the ideal is generated by the coefficients of the given polynomial in $t$. It is naturally an algebra over $H^*_{\GL_{2n}}\cong \C[c_1,\dots,c_{2n}]$. The action of $\theta$ is easy to figure out in this presentation, namely, all elements of degree $i$ will get multiplied by $(-1)^i$. To compute the fixed point scheme \eqref{fixed} we will have to determine  the coinvariant algebra of this $\theta$-action. To form the coinvariant algebra we add to the ideal the algebra elements which are acted upon by $\theta$ as $-1$, that is we add  the odd degree generators to the ideal. This way we get the following presentation of the coinvariant algebra: $H^{2*}_{\GL_{2n}}({\rm Gr}_{2k}(\C^{2n}),\C))_\theta\cong$
$$ \frac{\C[e_2,\dots,e_{2k}, f_2,\dots,f_{2n-2k}, c_2,\dots,c_{2n}]}{\left( (t^{2k}+e_2t^{2k-2}+\dots +e_{2k})(t^{2n-2k}+f_2t^{2n-2k-2}+\dots +f_{2n-2k})-(t^{2n}+c_2t^{2n-2}+\dots+c_{2n})\right)},$$  where all generators have even degree, and the ideal is generated by the coefficients of the indicated polynomial in $t^2$. This is an algebra over $\C[c_2,\dots,c_{2n}]$ which we can and will identify with $H^*_{{\rm Sp}(n)}:=H^*(B{\rm Sp}(n),\C)$ the cohomology ring of the classifying space of the compact unitary symplectic group. In turn, we can identify $H^{2*}_{\GL_{2n}}({\rm Gr}_{2k}(\C^{2n}),\C))_\theta$ with the equivariant cohomology $H^*_{{\rm Sp}(n)}(\Gr_k({\mathbb H}^n),\C)$ of the quaternionic Grassmannian $\Gr_k({\mathbb H}^n)$ of $k$-dimensional ${\mathbb \H}$-subspaces of ${\mathbb H}^n$, which is a non-Hermitian compact homogeneous  space isomorphic to $\Sp(n)/\Sp(k)\times \Sp(n-k)$. 

We can summarize our observation in the following diagram:

\begin{align}\label{thetagrassmann}\begin{array}{ccc} \Spec(H^{2*}_{\GL_{2n}}({\Gr}_{2k}(\C^{2n}),\C))^\theta &\cong & \Spec(H^{2*}_{\Sp(n)}(\Gr_{k}(\H^n),\C)) \\  \downarrow & & \downarrow \\ \Spec(H^{2*}_{\GL_{2n}})^\theta & \cong & \Spec(H^{2*}_{\Sp(n)})  \end{array} \end{align}
Thus, in light of \eqref{thetamodel}, the Hitchin map $h^\theta_{2k}$ on the even upward flow $W_{2k}^{2+}$ can be modelled by the spectrum of the equivariant cohomology of the quaternionic Grassmannian $\Gr_{k}(\H^n)$.

In fact, we can find similar coincidences of coinvariant algebras of the equivariant cohomology of cominuscule flag varieties in some other types. Cominuscule flag varieties correspond to maximal parabolic subgroups associated to minuscule coweights, or equivalently to simple roots, which occur with coefficient $1$ in the highest root (see also \cite[\S 8]{HauselHitchin} for more context).  For example we can consider the action of $\theta$ on \begin{align} \label{2*}H^{2*}_{\SO(4n+2)}(\SO(4n+2)/\SO(2)\times \SO(4n),\C)\end{align} given by $(-1)^{deg}$, where $deg=*$ is the degree of the grading on \eqref{2*}. The Hermitian symmetric space $\SO(4n+2)/\SO(2)\times \SO(4n)$ is an even quadric, a cominuscule flag variety for the special orthogonal group $\SO_{4n+2}$. The corresponding coinvariant algebra $$H^{2*}_{\SO(4n+2)}(\SO(4n+2)/\SO(2)\times \SO(4n),\C)_\theta\cong H^{2*}_{\SO(4n+1)}(\SO(4n+1)/\SO(4n),\C),$$ can be identified with the $\SO(4n+1)$-equivariant cohomology ring of the sphere $$S^{4n}\cong \SO(4n+1)/\SO(4n).$$ In fact we have the following diagram
\begin{align}\label{thetaquadric}\begin{array}{ccc} H^{2*}_{\SO(4n+2)}(\SO(4n+2)/\SO(2)\times \SO(4n),\C)_\theta&\cong &H^{2*}_{\SO(4n+1)}(\SO(4n+1)/\SO(4n),\C) \\ \uparrow & & \uparrow \\ \left(H^{2*}_{\SO(4n+2)}\right)_\theta & \cong & H^{2*}_{\SO(4n+1)}\end{array}.\end{align} In other words we can expect that the Hitchin map on even cominuscule upward flows in the $\SO_{4n+2}$ Higgs moduli space -- equivalently the Hitchin map on cominuscule upward flows in the $\SO(2n+2,2n)$-Higgs moduli space -- to be modelled by the spectrum of the equivariant cohomology of the sphere $S^{4n}$. 

For our final example we can consider the unique cominuscule flag variety for the exceptional ${\rm E}_6$. It is the complex Cayley plane ${\rm E}_6/{\rm Spin}({10})\times {\rm U}(1)$ which is  a compact Hermitian symmetric space. 
We can identify the coinvariant algebra of the $\theta=(-1)^{deg}$ action on its equivariant cohomology ring: 

\begin{align}\label{thetacayley}\begin{array}{ccc} H^{2*}_{{\rm E}_6}({\rm E}_6/{\rm Spin}({10})\times {\rm U}(1),\C)_\theta&\cong &H^{2*}_{{\rm F}_4}({\rm F}_4/{\rm Spin}(9),\C) \\ \uparrow & & \uparrow \\ \left(H^{2*}_{{\rm E}_6}\right)_\theta & \cong & H^{2*}_{{\rm F}_4}\end{array}.\end{align}

Thus we expect to model the Hitchin map on even cominuscule upward flows in the ${\rm E}_6$ Higgs moduli space -- equivalently on the cominuscule upward flows in the $E_{6(2)}$-Higgs moduli space -- by the spectrum of equivariant cohomology of the real Cayley plane ${\rm F}_4/{\rm Spin}(9)$. 
	
Mysteriously, in the above examples the symmetric spaces whose equivariant cohomology we found to give the $\theta$-coinvariant algebra of the equivariant cohomology of the cominuscule flag variety are homogeneous spaces for the Nadler group \cite[Table 1]{Nadler} of the corresponding quasi-split real form of Hodge type. That is $${\rm U}(n,n)^\vee \cong \Sp_n, \SO(2n+2,2n)^\vee \cong \SO_{4n+1}, \mbox{ and } {\rm E}_{6(2)}^\vee\cong {\rm F}_4.$$ Conjecturally, \cite[\S 7]{BaragliaSchaposnik} the Higgs bundle moduli space  for the Nadler group should give the support of the mirror of the Lagrangian brane given by the Higgs moduli space attached to a real form. However these appearances of the Nadler group remain to be understood. 

Even more surprising is that the Nadler groups in the above examples happen to be the fixed point subgroups of an involution of the ambient Langlands dual group. The corresponding anti-holomorphic involution can then be used to construct \cite{Elkner}  an anti-holomorphic involution in the cominuscule flag varieties in the above examples, and show that the Lagrangian fixed point manifold is isomorphic to the corresponding non-Hermitian compact symmetric spaces we have found above. Thus by uniformly denoting this anti-holomorphic involution by $\tau$ we get that $$\Gr_{2k}(\C^{2n})^\tau\cong \Gr_k(\H^n)$$ $$(\SO(4n+2)/\SO(2)\times \SO(4n))^\tau\cong \SO(4n+1)/\SO(4n)$$ and $$({\rm E}_6/{\rm Spin}({10})\times {\rm U}(1))^\tau \cong {\rm F}_4/{\rm Spin}(9).$$ Using Lemma~\ref{taustarcoinvariant} it can be shown - by observing that precisely the even degree generators of  invariant polynomials survive for the $\tau$-fixed groups - that in each of these cases $\tau$ induces our \begin{align} \label{taustar}
    \tau^*=\theta=(-1)^{deg} \end{align} on equivariant cohomology, where again $deg$ is half the degree of a homogeneous cohomology class. In turn, this observation and again Lemma~\ref{taustarcoinvariant} can be used to geometrically prove our final 

\begin{theorem}\label{final} The diagrams in \eqref{thetagrassmann}, \eqref{thetaquadric}, \eqref{thetacayley} commute and are induced by the involution $\tau$. 
\end{theorem}

\begin{remark} Because of \eqref{taustar} we can deduce, by the Lefschetz fixed point theorem, that the signatures of our Hermitian symmetric spaces $X=\G/\rm{H}$ agree \begin{align} \label{sign} {\rm sign}(X)=\tr(\theta:H^*(X)\to H^*(X))=\chi(X^\tau) \end{align} with the Euler characteristic of the corresponding non-Hermitian symmetric spaces $X^\tau$. The quantity ${\rm sign}(X)$ is relevant in our considerations as it agrees with the rank of the coinvariant algebra $H^*_\G(X)_\theta$, which in turn should compute the multiplicity of the even Hitchin map on the corresponding even cominuscule upward flow.

It is interesting to note that a very similar approach to \eqref{sign} was studied in \cite[Remark (1) p.337 ]{HirzebruchSlodowy} to determine the signature of our Hermitian symmetric spaces as the Euler characteristic of $X^\sigma$,  using the involution $\sigma:X\to X$ induced by the split real form - instead of our real form given by the Nadler group. The fixed point sets are different in the type $A$-case - certain real Grassmannians - from our quaternionic Grassmannians but the induced actions on the cohomology $\sigma^*=\tau^*=
\theta$ agree, because the real forms $\sigma$ and $\tau$ are inner to each other. 
\end{remark}
\begin{remark}
The involution $\sigma$ however suggests a solution for modelling the even Hitchin system inside the equivariant cohomology of the last family of cominuscule flag varieties with non-zero signature. Namely, we can consider the anti-holomorphic action of  $\sigma$ on  $\Gr_{2k}(\C^{2n+1})$ with fixed point set $\Gr_{2k}(\C^{2n+1})^\sigma\cong \Gr_{2k}(\R^{2n+1})$. Then it appears that we have $\sigma$ inducing \begin{align}\label{thetarealgrassmann} \begin{array}{ccc} H^{2*}_{\SL_{2n+1}}(\Gr_{2k}(\C^{2n+1}),\C)_\theta&\cong &H^{2*}_{\SO_{2n+1}}(\Gr_{2k}(\R^{2n+1}),\C) \\ \uparrow & & \uparrow \\ \left(H^{2*}_{\SL_{2n+1}}\right)_\theta & \cong & H^{2*}_{\SO_{2n+1}}\end{array}.\end{align} The subtlety of this case is that $\Gr_{2k}(\R^{2n+1})$ is no longer simply connected, thus the usual computation of its equivariant cohomology \eqref{eqcoh} does not apply. One can proceed by first determining the equivariant cohomology ring of its universal double cover - the oriented Grassmannian - from \eqref{eqcoh}, and then take invariants of the cover map. In fact, recently the  equivariant cohomology of $\Gr_{2k}(\R^{2n+1})$
was computed in \cite[Theorem 5.23]{He}, and the result matches  the coinvariant algebra $H^{2*}_{\SL_{2n+1}}(\Gr_{2k}(\C^{2n+1}),\C)_\theta$. We note that in this example the group $\SO_{2n+1}$ is the Langlands dual of the Nadler group $\Sp_n$ of our quasi-split real form $\SU(n,n+1)$ of Hodge type. In particular, their classifying spaces have isomorphic cohomology \begin{align} \label{long} H^*_{\SO_{2n+1}}\cong \C[\t_{\SO_{2n+1}}]^{{\rm W}_{\SO_{2n+1}}} \cong \C[\t_{\Sp_n}]^{{\rm W}_{\Sp_n}} \cong H^*_{\Sp_n},\end{align} because ${\rm W}_{\SO_{2n+1}}\cong {\rm W}_{\Sp_n}$, $\t_{\SO_{2n+1}}\cong \t_{\Sp_n}^*$ and the two representations in \eqref{long} can be identified by the Killing form. \end{remark}

\begin{remark} The final example of symmetric pairs from Remark~\ref{coinvariantsymmetric} we have not discussed yet is $(\SO_{4n},\SO_{4n-1})$. In this case  $\tau$ will induce on the cohomology of our cominuscule flag variety an involution $$\theta_\tau:=\tau^*:H^*_{\SO_{4n}}(\SO_{4n}/\SO_2\times \SO_{4n-2};\C)\to H^*_{\SO_{4n}}(\SO_{4n}/\SO_2\times \SO_{4n-2};\C),$$  which is different from the usual Hodge type $(-1)^{deg}$ unlike in the previous cases \eqref{taustar}. In fact we expect that this $\theta_\tau$ will be the involution corresponding to the (only) quasi-split real form $\SO(2n+1,2n-1)$ which is not split or of Hodge type.
Thus we expect that the spectrum of the diagram we get from Lemma~\ref{taustarcoinvariant} 
\begin{align}\label{thetatauquadric}\begin{array}{ccc} H^{2*}_{\SO_{4n}}(\SO_{4n}/\SO_2\times \SO_{4n-2};\C)_{\theta_\tau}&\cong &H^{2*}_{\SO_{4n-1}}(\SO_{4n-1}/\SO_{4n-2};\C) \\ \uparrow & & \uparrow \\ \left(H^{2*}_{\SO_{4n}}\right)_{\theta_\tau} & \cong & H^{2*}_{\SO_{4n-1}}\end{array}\end{align} models the Hitchin map on a cominuscule upward flow in the $\SO(2n+1,2n-1)$-Higgs moduli space.
    \end{remark}

\begin{remark}  Oscar Garc\'ia-Prada has pointed out to us that the symmetric pairs in Remark~\ref{coinvariantsymmetric}  precisely correspond to the complexifications of the maximal split subgroups of the non-split quasi-split real forms in \cite[Table 1. p. 2914]{Garcia-Prada-Peon-Nieto-Ramanan}. They are used \cite[Theorem 6.13.(1)]{Garcia-Prada-Peon-Nieto-Ramanan} to construct the Hitchin-Kostant-Rallis section in the  quasi-split real form cases. This demystifies their appearance in our descriptions of the Hitchin maps on cominuscule upward flows in the quasi-split Higgs moduli spaces in the diagrams \eqref{thetagrassmann}, \eqref{thetaquadric}, \eqref{thetacayley},\eqref{thetarealgrassmann} and \eqref{thetatauquadric} above. 
\end{remark}

\bibliographystyle{plain}
\bibliography{bib}

\begin{thebibliography}{10}

\bibitem{atiyah-bott}
Michael~F. Atiyah and Raoul Bott.
\newblock The moment map and equivariant cohomology.
\newblock {\em Topology. An International Journal of Mathematics}, 23(1):1--28,
  1984.

\bibitem{BaragliaSchaposnik}
David Baraglia and Laura~P. Schaposnik.
\newblock Real structures on moduli spaces of {H}iggs bundles.
\newblock {\em Advances in Theoretical and Mathematical Physics},
  20(3):525--551, 2016.

\bibitem{BialynickiBirula}
Andrzej Bialynicki-Birula.
\newblock Some {T}heorems on {A}ctions of {A}lgebraic {G}roups.
\newblock {\em Annals of Mathematics}, 98(3):480--497, 1973.

\bibitem{bradlow_surface_2003}
Steven~B. Bradlow, Oscar García-Prada, and Peter~B. Gothen.
\newblock Surface {Group} {Representations} and {U}(p, q)-{Higgs} {Bundles}.
\newblock {\em Journal of Differential Geometry}, 64(1), January 2003.

\bibitem{Elkner}
Mischa Elkner.
\newblock On pairs of symmetric spaces related by equivariant cohomology.
\newblock (I{S}{T}{A} rotation project 2023).

\bibitem{Franco-etal}
Emilio Franco, Peter~B. Gothen, André~G. Oliveira, and Ana Peón-Nieto.
\newblock {Narasimhan--Ramanan branes and wobbly Higgs bundles }.
\newblock arXiv:2302.02736.

\bibitem{HiggsPairs}
Oscar Garc\'ia-Prada, Peter~B. Gothen, and Ignasi Mundet~i Riera.
\newblock {The Hitchin-Kobayashi correspondence, Higgs pairs and surface group
  representations}.
\newblock 2009.
\newblock \href{https://arxiv.org/abs/0909.4487}{arXiv:0909.4487}.

\bibitem{Garcia-Prada-Peon-Nieto-Ramanan}
Oscar Garc\'{\i}a-Prada, Ana Pe\'{o}n-Nieto, and S.~Ramanan.
\newblock Higgs bundles for real groups and the {H}itchin-{K}ostant-{R}allis
  section.
\newblock {\em Trans. Amer. Math. Soc.}, 370(4):2907--2953, 2018.

\bibitem{garciaprada_involutions_2019}
Oscar García‐Prada and Sundararaman Ramanan.
\newblock Involutions and higher order automorphisms of {Higgs} bundle moduli
  spaces.
\newblock {\em Proceedings of the London Mathematical Society},
  119(3):681--732, September 2019.

\bibitem{goresky-kottwitz-macpherson}
Mark Goresky, Robert Kottwitz, and Robert MacPherson.
\newblock Equivariant cohomology, {K}oszul duality, and the localization
  theorem.
\newblock {\em Inventiones Mathematicae}, 131(1):25--83, 1998.

\bibitem{greub-etal}
Werner Greub, Stephen Halperin, and Ray Vanstone.
\newblock {\em Connections, curvature, and cohomology}.
\newblock Pure and Applied Mathematics, Vol. 47-III. Academic Press [Harcourt
  Brace Jovanovich, Publishers], New York-London, 1976.
\newblock Volume III: Cohomology of principal bundles and homogeneous spaces.

\bibitem{Hausel22}
Tam\'as Hausel.
\newblock Hitchin map on very stable and even very stable upward flows.
\newblock \url{https://www.youtube.com/watch?v=hWykyKDAWrU}, September 2022.

\bibitem{Hausel-hmsec}
Tam{\'a}s Hausel.
\newblock Hitchin map as spectrum of equivariant cohomology and {K}irillov
  algebras.
\newblock (in preparation).

\bibitem{HauselHitchin}
Tam{\'a}s Hausel and Nigel~J. Hitchin.
\newblock Very stable {H}iggs bundles, equivariant multiplicity and mirror
  symmetry.
\newblock {\em Inventiones mathematicae}, 228(2):893--989, May 2022.

\bibitem{He}
Chen He.
\newblock Localization of equivariant cohomology rings of real {G}rassmannians.
\newblock arXiv:1609.06243.

\bibitem{HirzebruchSlodowy}
Friedrich Hirzebruch and Peter Slodowy.
\newblock Elliptic genera, involutions, and homogeneous spin manifolds.
\newblock {\em Geometriae Dedicata}, 35(1-3):309--343, 1990.

\bibitem{HitchinIntegrable}
Nigel~J. Hitchin.
\newblock {Stable bundles and integrable systems}.
\newblock {\em Duke Mathematical Journal}, 54(1):91 -- 114, 1987.

\bibitem{HitchinSelfDuality}
Nigel~J. Hitchin.
\newblock {The Self-Duality Equations on a Riemann Surface}.
\newblock {\em Proceedings of the London Mathematical Society},
  s3-55(1):59--126, 1987.

\bibitem{HitchinSection}
Nigel~J. Hitchin.
\newblock Lie groups and {{T}}eichmüller space.
\newblock {\em Topology}, 31(3):449--473, 1992.

\bibitem{HitchinMultiplicity}
Nigel~J. Hitchin.
\newblock Multiplicity algebras for rank 2 bundles on curves of small genus.
\newblock arXiv:2203.03424.

\bibitem{Laumon}
G{\'e}rard Laumon.
\newblock {Un analogue global du cône nilpotent}.
\newblock {\em Duke Mathematical Journal}, 57(2):647 -- 671, 1988.

\bibitem{Nadler}
David Nadler.
\newblock Perverse sheaves on real loop {G}rassmannians.
\newblock {\em Inventiones Mathematicae}, 159(1):1--73, 2005.

\bibitem{Nitsure}
Nitin Nitsure.
\newblock {Moduli Space of Semistable Pairs on a Curve}.
\newblock {\em Proceedings of the London Mathematical Society},
  s3-62(2):275--300, 1991.

\bibitem{PanyushevYakimova}
Dmitri~I. Panyushev and Oksana~S. Yakimova.
\newblock Automorphisms of finite order, periodic contractions, and
  {P}oisson-commutative subalgebras of {$ S(\mathfrak{g})$}.
\newblock {\em Math. Z.}, 303(2):Paper No. 51, 30, 2023.

\bibitem{Peon-Nieto}
Ana Peón-Nieto.
\newblock Wobbly moduli of chains, equivariant multiplicity and ${\U}(n_0,
  n_1)$-{H}iggs bundles.

\bibitem{BNR}
Sundararaman Ramanan, Mudumbai~S. Narasimhan, and Arnaud Beauville.
\newblock {Spectral curves and the generalised theta divisor.}
\newblock {\em Journal für die reine und angewandte Mathematik}, 398:169--179,
  1989.

\bibitem{Simpson}
Carlos~T. Simpson.
\newblock {Moduli of representations of the fundamental group of a smooth
  projective variety I}.
\newblock {\em Publications Math{\'e}matiques de l'Institut des Hautes
  {\'E}tudes Scientifiques}, 79(1):47--129, 1994.

\bibitem{SimpsonII}
Carlos~T. Simpson.
\newblock Moduli of representations of the fundamental group of a smooth
  projective variety {II}.
\newblock {\em Publications Math{\'e}matiques de l'Institut des Hautes
  {\'E}tudes Scientifiques}, 80:5--79, 1994.

\bibitem{Springer}
T.~A. Springer.
\newblock Regular elements of finite reflection groups.
\newblock {\em Invent. Math.}, 25:159--198, 1974.

\end{thebibliography}
\end{document}